\documentclass[11pt]{article}
\usepackage[utf8]{inputenc}
\usepackage{amssymb}%
\usepackage{amsthm}%
\usepackage{amsmath}
\usepackage{amsfonts}
\usepackage{dsfont}
\usepackage{stmaryrd}
\usepackage[table,xcdraw]{xcolor}%
\usepackage[labelformat=simple]{subfig}

\usepackage[font=footnotesize]{caption}
\usepackage{graphicx}
\usepackage{authblk}
\usepackage{natbib}
\usepackage{booktabs}
\usepackage{multirow}

\usepackage{changepage}

\usepackage{hyperref}
\hypersetup{%
        colorlinks,
        breaklinks=true,
        plainpages=false,%
        citecolor=blue,
        linkcolor=blue,
        urlcolor=blue,
        bookmarksopen=true,%
        bookmarksnumbered=false,%
        bookmarksdepth=5%
}

\usepackage{xurl}

\usepackage[shortlabels]{enumitem}

\usepackage[export]{adjustbox}
\usepackage{fullpage}

\newcommand{\norm}[1]{\left\lVert#1\right\rVert}
\newcommand{\parent}[1]{\left(#1\right)}

\newcommand{\dotprod}[1]{\left<#1\right>}
\renewcommand{\brack}[1]{\left[#1\right]}
\renewcommand{\brace}[1]{\left\{#1\right\}}
\renewcommand{\det}[1]{\left|#1\right|}
\newcommand{\floor}[1]{\left\lfloor#1\right\rfloor}

\newcommand{\argmax}[1]{\underset{#1}{\textnormal{argmax}} \,}

\renewcommand{\min}[1]{\underset{#1}{\textnormal{min}} \,}
\renewcommand{\max}[1]{\underset{#1}{\textnormal{max}} \,}
\renewcommand{\sup}[1]{\underset{#1}{\textnormal{sup}} \,}

\newcommand{\R}{\mathbb{R}}

\newcommand{\N}{\mathbb{N}}
\renewcommand{\S}{\mathcal{S}}
\renewcommand{\L}{\mathcal{L}}

\providecommand{\keywords}[1]
{
  \small	
  \textbf{\textit{Keywords---}} #1
}

\newtheorem{Theorem}{Theorem}
\newtheorem{Proposition}[Theorem]{Proposition}%
\newtheorem{Lemma}[Theorem]{Lemma}

\newtheorem{Remark}{Remark}%

\begin{document}
\title{Deterministic Approximate EM Algorithm; Application to the Riemann Approximation EM and the Tempered EM}

\author[1, 2]{Thomas Lartigue}
\author[1]{Stanley Durrleman}
\author[3]{St\'ephanie Allassonni\`ere}
\affil[1]{Aramis Project-Team, Inria, 75012 Paris, France}
\affil[2]{CMAP, CNRS, \'Ecole Polytechnique, Institut Polytechnique de Paris, 91120 Palaiseau, France}
\affil[3]{HeKA, Centre de recherche des Cordeliers, Universit\'e de Paris, INRIA, INSERM, Sorbonne Universit\'e, 75012 Paris, France}
\date{}                     
\setcounter{Maxaffil}{0}
\renewcommand\Affilfont{\itshape\small}

\maketitle

\begin{abstract}
    The Expectation Maximisation (EM) algorithm is widely used to optimise non-convex likelihood functions with latent variables. Many authors modified its simple design to fit more specific situations. For instance, the Expectation (E) step has been replaced by Monte Carlo (MC), Markov Chain Monte Carlo or tempered approximations, etc. Most of the well-studied approximations belong to the stochastic class. By comparison, the literature is lacking when it comes to deterministic approximations. In this paper, we introduce a theoretical framework, with state-of-the-art convergence guarantees, for any deterministic approximation of the E step. We analyse theoretically and empirically several approximations that fit into this framework. First, for intractable E-steps, we introduce a deterministic version of MC-EM using Riemann sums. A straightforward method, not requiring any hyper-parameter fine-tuning, useful when the low dimensionality does not warrant a MC-EM. Then, we consider the tempered approximation, borrowed from the Simulated Annealing literature and used to escape local extrema. We prove that the tempered EM verifies the convergence guarantees for a wider range of temperature profiles than previously considered. We showcase empirically how new non-trivial profiles can more successfully escape adversarial initialisations. Finally, we combine the Riemann and tempered approximations into a method that accomplishes both their purposes.
\end{abstract}



\section{Introduction}
The Expectation Maximisation (EM) algorithm was introduced by Dempster, Laird and Rubin (DLR, \cite{dempster1977maximum}) to maximise likelihood functions $g(\theta)$ defined from inherent latent variables $z$ and that were non-convex and had intricate gradients and Hessians. The EM algorithm estimates $\theta$ in an iterative fashion, starting from a certain initial value $\theta_0$ and where the new estimate $\theta_{n+1}$ at step $n+1$ is function the estimate $\theta_n$ from the previous step $n$. In addition to presenting the method, DLR provide convergence guarantees on the sequence of estimated parameters $\{\theta_n\}_n$, namely that it converges towards a critical point of the likelihood function as the step $n$ of the procedure grows. Their flawed proof for this result was later corrected by \citet{wu1983convergence}, and more convergence guarantees were studied by \citet{boyles1983convergence}. 
Since some likelihood functions are too complex to apply DLR's raw version of the EM, authors of later works have proposed alternative versions, usually with new convergence guarantees. On the one hand, when the maximisation step (M step) is problematic, other optimisation methods such as coordinate descent \cite{wu1983convergence} or gradient descent \cite{lange1995gradient} have been proposed.
On the other hand, several works introduce new versions of the algorithm where the expectation step (E step), which can also be intractable, is approximated. Most of them rely on Monte Carlo (MC) methods and Stochastic Approximations (SA) to estimate this expectation. Notable examples include the Stochastic Approximation EM (SAEM, \cite{delyon1999convergence}), the Monte Carlo EM (MC-EM, \cite{wei1990monte}), the Markov Chain Monte Carlo EM (MCMC-EM, \cite{fort2003convergence}), the MCMC-SAEM \cite{kuhn2005maximum, allassonniere2010construction} and the Approximate SAEM \cite{allassonniere2021new}.
Random approximation of the E step have also been used in the case where the data are too voluminous to allow a full E step computation. {A non-exhaustive list includes} the Incremental EM \cite{neal1998view, ng2003choice}, the Online EM \cite{cappe2009line} and, more recently, the stochastic EM with variance reduction (sEM-vr) \cite{chen2018stochastic}, the fast Incremental EM (FIEM) \cite{karimi2019global}, the Stochastic Path-Integrated Differential EstimatoR EM (SPIDER-EM) \cite{fort2020stochastic}, {and the mini-batch MCMC-SAEM} \cite{kuhn2020properties}. Most of these variants come with their own theoretical convergence guarantees for the models of the exponential family. Recent works have also provided theoretical analysis of the EM algorithm outside of the exponential family, with locally strongly-concave log-likelihood function around the global maxima by \citet{balakrishnan2017statistical}, or without such strong-concavity assumption by \citet{dwivedi2020singularity}.

The stochastically approximated EM algorithms constitute an extensive catalogue of methods. Indeed, there are many possible variants of MCMC samplers that can be considered, as well as the additional parameters, such as the burn-in period length and the gain decrease sequence, that have to be set. All these choices have an impact on the convergence of these ``EM-like'' algorithms and making the appropriate ones for each problem can be overwhelming, see \cite{booth1999maximizing, levine2001implementations, levine2004automated}, among others, for discussions on tuning the MC-EM alone. On several cases, one might desire to dispose of a simpler method, possibly non-stochastic, and non-parametric to run an EM-like algorithm for models with no closed forms. However the literature is lacking in that regards. The Quasi-Monte Carlo EM, introduced by \citet{pan1998quasi}, is a deterministic version of Monte Carlo EM, however theoretical guarantees are not provided. In that vein, \citet{jank2005quasi} introduces the randomised Quasi-Monte Carlo EM, which is not deterministic, and does not have theoretical guarantees either. {Another example of deterministic approximation can be found in the Variational \mbox{EM \cite{attias1999Inferring, bishop2006pattern, tzikas2008variational}}, where the conditional distribution is approximated within a certain distribution family, and optimised over before each E step. Note that with such a design, if the distance between the true conditional and the chosen family is non-zero, then the approximation can never converge towards the true conditional.}
 
In addition to intractable E steps, EM procedures also commonly face a second issue: their convergence, when guaranteed, can be towards any maximum. This theoretical remark has crucial numerical implications. Indeed, most of the time, convergence is reached towards a sub-optimal local maximum, usually very dependent on the initialisation. To tackle this issue and improve the solutions of the algorithm, other types of, usually deterministic, approximations of the E step have been proposed. One notable example is the tempering (or ``annealing'') of the conditional probability function used in the E step. Instead of replacing an intractable problem by a tractable one, the tempering approximation is used to find better local maxima of the likelihood profile during the optimisation process, in the spirit of the Simulated Annealing of \citet{kirkpatrick1983optimization} and the Parallel Tempering (or Annealing MCMC) of \citet{swendsen1986replica, geyer1995annealing}. The deterministic annealing EM was introduced by \citet{ueda1998deterministic} with a decreasing temperature profile; another temperature profile was proposed by \mbox{\citet{naim2012convergence}}. Contrary to most of the studies on stochastic approximations, these two works do not provide theoretical convergence guarantees for the proposed tempered methods, which, as a consequence, does not provide insight on the choice of the temperature scheme. Moreover, the tempered methods do not allow the use of the EM in case of an intractable E step. In their tempered SAEM algorithm, \citet{allassonniere2021new} combine the stochastic and tempering approximations, which allows the SAEM to run, even with an intractable E step, while benefiting from the improved optimisation properties brought by the tempering. In addition, theoretical convergence guarantees are provided. However, this method is once again stochastic and parametric.

Overall, most of the literature on approximated E steps focuses on stochastic approximations that estimate intractable conditional probability functions. The few purely deterministic approximations proposed, such as the tempered/annealed EM, are used for other purposes, improving the optimisation procedure, and lack convergence guarantees.

{
In this paper, we prove that, with mild model assumptions and regularity conditions on the approximation, any Deterministic Approximate EM benefits from the state of the art theoretical convergence guarantees of \cite{wu1983convergence, lange1995gradient, delyon1999convergence}. This theorem covers several already existing methods, such as the tempered EM, and paves the way for the exploration of new ideas. To illustrate this, we study a small selection of methods with practical applications that verify our theorem's hypotheses. 
}
First, for E steps without closed form, we propose to use Riemann sums to estimate the intractable normalising factor. This ``Riemann approximation EM'' is a deterministic, less parametric, alternative to the MC-EM and its variants. {It is suited to the small dimension cases, where the analytic estimation of the intractable integral is manageable, and MC estimation unnecessary.} Second, we prove that the deterministic annealed EM (or ``tempered EM'') of \citet{ueda1998deterministic} is a member of our general deterministic class as well. We prove that the convergence guarantees are achieved with almost no condition of the temperature scheme, justifying the use of a wider range of temperature profile than those proposed by \cite{ueda1998deterministic} and \cite{naim2012convergence}. Finally, since the Riemann and tempered approximations are two separate methods that fulfil very different practical purposes, we also propose to associate the two approximations in the ``tempered Riemann approximation EM'' when both their benefits are desired.

In Section \ref{sect:approx_EM}, we introduce our general class of deterministic approximated versions of the EM algorithm and prove their convergence guarantees, for models of the exponential family. We discuss the ``Riemann approximation EM'' in Section \ref{sect:riemann_em}, the ``tempered EM'' in Section \ref{sect:tmp_em}, and their association, ``tempered Riemann approximation EM'', in \mbox{Section \ref{sect:tmp_riemann_em}}. In Section \ref{sect:exp}, we study empirically each of these three methods. Section \ref{sect:conclusions} discusses and concludes this work.

We demonstrate empirically that the Riemann EM converges properly on a model with and an intractable E step, and that adding the tempering to the Riemann approximation allows in addition to get away from the initialisation and recover the true parameters. On a tractable Gaussian Mixture Model, we compare the behaviours and performances of the tempered EM and the regular EM. In particular, we illustrate that the tempered EM is able to escape adversarial initialisations, and consistently reaches better values of the likelihood than the unmodified EM, in addition to better estimating the model parameters.


\section{Materials and Methods}

\subsection{Deterministic Approximate EM Algorithm and Its Convergence for the Exponential Family} \label{sect:approx_EM}
\subsubsection{Context and Motivation}\label{sect:EM_context}
In this section, we propose a new class of deterministic EM algorithms with approximated E step. {This class of algorithms is very general. In particular, it includes methods with deterministic approximations of intractable normalisation constant. It also includes optimisation-oriented approximations, such as the annealing approximation used to escape sub-optimal local extrema of the objective. In addition, combinations of these two types of approximations are also covered.} We prove that members of this class benefit from the same convergence guarantees that can be found in the state of the art references \citep{wu1983convergence, lange1995gradient, delyon1999convergence} for the classical EM algorithm, and under similar model assumptions. The only condition on the approximated distribution being that it converges towards the real conditional probability distribution with a $l_2$ regularity. Like the authors of \cite{delyon1999convergence, fort2003convergence, allassonniere2021new}, we work with probability density functions belonging to the exponential family. The specific properties of which are providen in the hypothesis $M1$ of Theorem \ref{thm:main}.

{EM algorithms are most often considered in the context of the \textit{missing data problem}, see \cite{delyon1999convergence}. The formulation of the problem is the following: we observe a random variable $x$, described by a parametric probability density function (pdf) noted $g(\theta)$ with parameter $\theta \in \Theta \subset \R^l, \; l\in\N^*$. We assume that there exists a hidden variable $z$ informing the behaviour of the observed variable $x$ such that $g(\theta)$ can be expressed as the integral of the complete likelihood $h(z; \theta)$: $g(\theta) = \int_z h(z; \theta) \mu(dz)$, with $\mu$ the reference measure. We denote $p_{\theta}(z) := h(z; \theta)/g(\theta)$ the conditional density function of $z$. For a given sample $x$, the goal is to maximise the likelihood $g(\theta)$ with respect to $\theta$. In all these notations, and throughout the rest of the paper, we omit $x$ as a variable since it is considered fixed once and for all, and everything is done conditionally to $x$. In particular, in many applications, the observed data $x$ is made of $N\in \N^*$ samples: $x := (x^{(1)}, x^{(N)})$. In such cases, the sample size $N < +\infty$ is supposed finite and fixed once and for all. We are never in the asymptotic regime $N \longrightarrow +\infty$.}

The foundation of the EM algorithm is that while $\ln g(\theta)$ is hard to maximise in $\theta$, the functions $\theta \mapsto \ln h(z; \theta)$ and even $\theta \mapsto \mathbb{E}_{z} \brack{\ln h(z; \theta)}$ are easier to work with because of the information added by the latent variable $z$ (or its distribution). In practice however, the actual value of $z$ is unknown and its distribution $p_{\theta}(z)$ dependent on $\theta$. Hence, the EM was introduced by DLR \cite{dempster1977maximum} as the two-stages procedure starting from an initial point $\theta_0$ and iterated over the number of steps $n$:
\begin{itemize}
 \item[] (E) 
With the current parameter $\theta_n$, calculate the conditional probability $p_{\theta_n}(z)$;
 \item[] (M) 
 To get $\theta_{n+1}$, maximise in $\theta \in \Theta$ the function $\theta \mapsto \mathbb{E}_{z \sim p_{\theta_n}(z)} \brack{\ln h(z; \theta)}$;
\end{itemize}

Which can be summarised as:
\begin{equation} \label{eq:EM_steps}
 \theta_{n+1} := T(\theta_n) := \argmax{\theta \in \Theta} \mathbb{E}_{z \sim p_{\theta_n}(z)} \brack{\ln h(z; \theta)} \, ,
\end{equation}
where we call $T$ the point to point map in $\Theta$ corresponding to one EM step.
We will not redo the basic theory of the exact EM here, but this procedure noticeably increase $g(\theta_n)$ at each new step $n$. However, as discussed in the introduction, in many cases, one may prefer to or need to use an approximation of $p_{\theta_n}(z)$ instead of the exact analytical value.

In the following, we consider a deterministic approximation of $p_{\theta}(z)$ noted $\tilde{p}_{\theta, n}(z)$ which depends on the current step $n$ and on which we make no assumption at the moment. The resulting steps, defining the ``Approximate EM'', can be written under the same form as \eqref{eq:EM_steps}:
\begin{equation} \label{eq:approximate_EM_steps}
 \theta_{n+1} := F_n(\theta_n) := \argmax{\theta \in \Theta} \mathbb{E}_{z \sim \tilde{p}_{\theta_n, n}(z)} \brack{\ln h(z; \theta)} \, ,
\end{equation}
where $\brace{F_n}_{n \in \N}$ is the sequence of point to point maps in $\Theta$ associated with the sequence of approximations $\brace{\tilde{p}_{\theta, n}(z)}_{ \theta \in \Theta; \, n \in \N}$. In order to ensure the desired convergence guarantees, we add a slight modification to this EM sequence: \textit{re-initialisation of the EM sequence onto increasing compact sets}. This modification was introduced by \citet{chen1987convergence} and adapted by \citet{fort2003convergence}. Assume that you dispose of an increasing sequence of compacts $\brace{K_n}_{n \in \mathbb{N}}$ such that $\cup_{n \in \mathbb{N}} K_n = \Theta$ and $\theta_0 \in K_0$. Define $j_0 := 0$. Then, the transition $\theta_{n+1} = F_n(\theta_n)$ is accepted only if $F_n(\theta_n)$ belongs to the current compact $K_{j_n}$, otherwise the sequence is reinitialised at $\theta_0$. The steps of the resulting algorithm, called \textit{Stable Approximate EM}, can be written as:
\begin{equation} \label{eq:stable_approximate_EM_steps}
 \begin{cases}
 if F_n(\theta_n)\! \in\! K_{j_n}, \, \text{then} \; \theta_{n+1} = F_n(\theta_n),\, \text{and} \; j_{n+1} := j_n\\
 if F_n(\theta_n)\! \notin\! K_{j_n}, \, \text{then} \; \theta_{n+1} = \theta_0,\, \text{and} \; j_{n+1} := j_n+1 \, .
 \end{cases}
\end{equation}

This re-initialisation of the EM sequence may seem like a hurdle, however, this truncation is {only} a theoretical requirement. In practice, the first compact $K_0$ is taken so large that it covers the most probable areas of $\Theta$ and the algorithms \eqref{eq:approximate_EM_steps} and \eqref{eq:stable_approximate_EM_steps} are identical as long as the sequence $\brace{\theta_n}_n$ does not diverges towards the border of $\Theta$.

\subsubsection{Theorem}
In the following, we will state the convergence Theorem of Equation~\eqref{eq:stable_approximate_EM_steps} and provide a brief description of the main steps of the proof.
\begin{Theorem}[Convergence of the Stable Approximate EM]\label{thm:main}
 Let $\brace{\theta_n}_{n \in \N}$ be a sequence of the Stable Approximate EM defined in Equation~\eqref{eq:stable_approximate_EM_steps}. Let us assume two sets of hypotheses:
 \begin{itemize}
 \item \textbf{Conditions M1--3 on the model.}
 \begin{itemize}
 \item[] \textbf{M1.} $\Theta \subseteq \R^l$, $\mathcal{X} \subseteq \R^d$ and $\mu$ is a $\sigma$-finite positive Borel measure on $\mathcal{X}$. Let $\psi :\Theta \to \R$, $\phi :\Theta \to \R^q$ and $S : \mathcal{X} \to \mathcal{S} \subseteq \R^q$. Define $L: \mathcal{S} \times \Theta \to \R$, $h : \mathcal{X} \times \Theta \to \R_+^*$ and $g :\Theta \to \R_+^*$ as:
 \begin{equation*}
 \begin{split}
 L(s; \theta ) &:= \psi(\theta) + \dotprod{s,\phi(\theta)} \, ,\\
 h(z; \theta ) &:= \exp( L(S(z); \theta ))\, , \\
 g(\theta) &:= \int_{z \in \mathcal{X}} h(z; \theta) \mu(dz) \, .
 \end{split}
 \end{equation*}
 \item[] \textbf{M2.} Assume that
 \begin{itemize}
 \item[] (a*) $\psi$ and $\phi$ are continuous on $\Theta$;
 \item[] (b) for all $\theta \in \Theta$, $\bar{S}(\theta) := \int_z S(z) p_{\theta}(z)\mu(dz)$ is finite and continuous on $\Theta$;
 \item[] (c) there exists a continuous function $\hat{\theta} :\mathcal{S} \to \Theta$ such that for all $s \in \S, \; L(s; \hat{\theta}(s)) = \sup{\theta\in \Theta} L(s; \theta)$;
 \item[] (d) $g$ is positive, finite and continuous on $\Theta$ and the level set $\brace{\theta \in \Theta, g(\theta) \geq M}$ is compact for any $M > 0$.
 \end{itemize}
 \item[] \textbf{M3.} {Define $\L \!:=\! \brace{\theta \!\in\! \Theta | \hat{\theta}\circ\bar{S}(\theta) = \theta}$ and, for any $g^*\in \R_+^*$, $\L_{g^*}\! :=\! \brace{\theta \!\in\! \L | g(\theta) \!= \!g^* }$}. Assume either that:
 \begin{itemize} 
 \item[] ($a$) The set $g(\L)$ is compact or
 \item[] ($a'$) for all compact sets $K \subseteq \Theta,$ $g\parent{K \cap \L}$ is finite.
 \end{itemize}

 \end{itemize}
 \item \textbf{The conditions on the approximation.} Assume that $\Tilde{p}_{\theta, n}(z)$ is deterministic. Let $S(z)= \brace{S_u(z)}_{u=1}^q$. For all indices $u$, for any compact set $K \subseteq \Theta$, one of the two following configurations holds:
 \begin{equation} \label{eq:first_sufficient_condition}
 \begin{cases}
 \int_z S^2_u(z) dz < \infty \, ,\\
 \sup{\theta \in K}\, \int_z \parent{ \Tilde{p}_{\theta, n}(z) - p_{\theta}(z) }^2 dz \underset{n \to \infty}{\longrightarrow} 0 \, .
 \end{cases}
 \end{equation}
 Or
 \begin{equation} \label{eq:second_sufficient_condition}
 \begin{cases}
 \sup{\theta \in K}\, \int_z S^2_u(z) p_{\theta}(z) dz < \infty \, ,\\
 \sup{\theta \in K}\, \int_z \parent{ \frac{\Tilde{p}_{\theta, n}(z)}{p_{\theta}(z)} - 1 }^2 p_{\theta}(z) dz \underset{n \to \infty}{\longrightarrow} 0 \, .
 \end{cases}
 \end{equation}
 \end{itemize}
Then, 
 \begin{enumerate}[(i)]
 \item \begin{enumerate}[(a)]
 \item $\underset{n \to \infty}{\lim}\, j_n < \infty$ and $ \sup{n \in \mathbb{N}}\norm{\theta_n} < \infty $, with probability 1;
 \item $g(\theta_n)$ converges towards a connected component of $g(\L)$.
 \end{enumerate}
 \item {Let $Cl : \mathcal{P}(\Theta) \to \Theta$ be the set closure function and $d : \Theta \times \mathcal{P}(\Theta) \to \R_+$ be any point-to-set distance function within $\Theta$.} If $g\parent{ \L \cap Cl\parent{\brace{\theta_n}_{n \in \mathbb{N}}} }$ has an empty interior, then $\exists g^* \in \R_+^*$ such that: 
 \begin{equation*}
 \begin{split}
 g(\theta_n) \underset{n \to \infty}{\longrightarrow} g^*\, ,\\
 d(\theta_n, \L_{g^*}) \underset{n \to \infty}{\longrightarrow} 0\, .
 \end{split}
 \end{equation*}
 \end{enumerate}
\end{Theorem}

\begin{Remark}
~~~

{\begin{itemize}
 \item {When $g$ is differentiable on $\Theta$, its stationary points are the stable points of the EM, i.e. $\L = \brace{\theta \in \Theta | \nabla g (\theta) = 0 }$.} 
 \item The M1--3 conditions in Theorem \ref{thm:main} are almost identical to the similarly named M1--3 conditions in \cite{fort2003convergence}. The only difference is in $M2$ $(a^*)$, where we remove the hypothesis that $S$ has to be a continuous function of $z$, since it is not needed when the approximation is not stochastic. 
 \item The property $\int_z S^2_u(z) dz < \infty$ of the condition \eqref{eq:first_sufficient_condition} can seem hard to verify since $S$ is not integrated here against a probability density function. However, when $z$ is a finite variable, as is the case for finite mixtures, this integral becomes a finite sum. Hence, condition \eqref{eq:first_sufficient_condition} is better adapted to finite mixtures, while condition \eqref{eq:second_sufficient_condition} is better adapted to continuous hidden variables.
 \item The two sufficient conditions \eqref{eq:first_sufficient_condition} and \eqref{eq:second_sufficient_condition} involve a certain form of $l_2$ convergence of $\tilde{p}_{\theta, n}$ towards $p_{\theta}$. If the latent variable $z$ is continuous, this excludes countable (and finite) approximations such as sums of Dirac functions, since they have a measure of zero. In particular, this excludes Quasi-Monte Carlo approximations. However, one look at the proof of Theorem \ref{thm:main} (in Appendix \ref{sect:proofs_main}) or at the following sketch of proof reveals it is actually sufficient to verify $\sup{\theta \in K}\, \norm{\Tilde{S}_n(\theta) - \bar{S}(\theta)} \underset{n \to \infty}{\longrightarrow} 0$ for any compact set $K$. {Where $\Tilde{S}_n(\theta) := \int_z S(z) \Tilde{p}_{\theta, n}(z)\mu(dz)$ denotes the approximated E step in the Stable Approximate EM. This condition can be verified by finite approximations.}
\end{itemize}}
\end{Remark}

\subsubsection{Sketch of Proof} {The detailed proof of this result can be found in Appendix \ref{sect:proofs_main}, we propose here an abbreviated version where we highlight the key steps.

Two intermediary propositions, introduced and proven in \citet{fort2003convergence}, are instrumental in the proof of Theorem \ref{thm:main}. These two propositions are called Propositions 9 and 11 by Fort and Moulines, and used to prove their Theorem 3. Within our framework, with the new condition $M2(a^*)$ and the absence of Monte Carlo sum, the reasoning for verifying the conditions of applicability of the two proposition is quite different from \cite{fort2003convergence} and will be highlighted below. Let us state these two propositions using the notations of this paper:
\begin{Proposition}[``Proposition 9'']\label{prop:9}{
 Consider a parameter space $\Theta \subseteq \R^l$, a point to point map $T$ on $\Theta$, a compact $K \subset \Theta$ and a subset $\L \subseteq \Theta$ such that $ \L \cap K $ is compact. Let $g$ be a $C^0$, Lyapunov function relative to $(T, \L)$. Assume that there exist a $K-$valued sequence $\brace{\theta_n}_{n \in \N^*} \in K^{\mathbb{N}}$ such that:
 \begin{equation*}
 \underset{n \to \infty}{\lim} |g(\theta_{n+1}) - g \circ T(\theta_n) | = 0\, .
 \end{equation*}
 Then 
 \begin{itemize}
 \item $\brace{g(\theta_n)}_{n \in \mathbb{N}} $ converges towards a connected component of $g(\L \cap K)$
 \item If $g(\L \cap K)$ has an empty interior, then $\exists g^* \in \R_+^*$ such that $\brace{g(\theta_n)}_n$ converges towards $g^*$. Moreover, $\brace{\theta_n}_n$ converges towards the set $\L_{g^*} \cap K$. Where $\L_{g^*} := \brace{\theta \in \L | g(\theta) = g^* }$.
 \end{itemize}}
\end{Proposition}
\begin{Proposition}[``Proposition 11'']\label{prop:11}{
 Consider a parameter space $\Theta \subseteq \R^l$, a subset $\L \subseteq \Theta$, a point to point map $T$ on $\Theta$ and a sequence of point to point maps $\brace{F_n}_{n\in \N^*}$ also on $\Theta$. Let $\brace{\theta_n}_{n\in \N} \in \Theta^{\N}$ be a sequence defined from $\brace{F_n}_n$ by the Stable Approximate EM equation \eqref{eq:stable_approximate_EM_steps} for some increasing sequence $\brace{K_n}_{n\in\N}$ of compacts of $\Theta$. Let $\brace{j_n}_{n\in \N}$ be the corresponding sequence of indices, also defined in \eqref{eq:stable_approximate_EM_steps}, such that $\forall n, \; \theta_n \in K_{j_n}$. We assume:
 \begin{itemize}
 \item[(a)] the $C1-2$ conditions of Proposition 10 of \cite{fort2003convergence}.
 \begin{itemize}
 \item \textbf{(C1)} There exists $g$, a $C^0$ Lyapunov function relative to $(T, \L)$ such that for all $M > 0$, $\brace{\theta \in \Theta | g(\theta) > M}$ is compact, and:
 \begin{equation*}
 \Theta = \cup_{n \in \N} \brace{\theta \in \Theta | g(\theta) > n^{-1}}\, .
 \end{equation*}
 \item \textbf{(C2)} $g(\L)$ is compact OR \textbf{(C2')} $g(\L \cap K)$ is finite for all compact $K \subseteq \Theta$.
 \end{itemize}
 \item[(b)] $\forall \theta \in K_0, \quad \underset{n \to \infty}{\lim}\,|g \circ F_n - g \circ T|(\theta) = 0$
 \item[(c)] $\forall$ compact $K \subseteq \Theta$: 
 \begin{equation*}
 \underset{n \to \infty}{\lim} |g \circ F_n(\theta_n) - g \circ T (\theta_n)| \mathds{1}_{\theta_n \in K} = 0\, .
 \end{equation*}
 \end{itemize}
 Then:

 With probability 1, $\underset{n \to \infty}{\textnormal{lim sup}}\, j_n < \infty$ and $\brace{\theta_n}_n$ is a compact sequence.
 }
\end{Proposition}
\begin{Remark}
 In \cite{fort2003convergence}, condition $C1$ of Proposition \ref{prop:11} is mistakenly written as: 
 $$\Theta = \cup_{n \in \N} \brace{\theta \in \Theta | W(\theta) > n}\, .$$
 This is a typo that we have corrected here.
\end{Remark}
The proof of Theorem \ref{thm:main} is structured as follows:
verifying the conditions of Proposition~\ref{prop:11}, applying Proposition \ref{prop:11}, verifying the conditions of Proposition \ref{prop:9} and finally applying Proposition \ref{prop:9}.

\vspace{12pt}

\noindent \textbf{Verifying the conditions of Proposition \ref{prop:11}.}
Let $g$ be the observed likelihood function defined in hypothesis $M1$ of Theorem \ref{thm:main}, $T$ the exact EM point to point map defined in \eqref{eq:EM_steps}, $\L := \brace{\theta \in \Theta | T(\theta) = \theta }$ the set of its stable points, $\brace{F_n}_n$ a sequence of approximated point to point map as defined in \eqref{eq:approximate_EM_steps}. With some increasing sequence $\brace{K_n}_{n\in\N}$ of compacts of $\Theta$, let $\brace{\theta_n, j_n}_{n\in \N}$ be defined from the Stable Approximate EM equation \eqref{eq:stable_approximate_EM_steps}.

By design of the regular EM, the following two properties are true. First, $g$ is a $C^0$, Lyapunov function relative to $(T, \L)$. Second, the map $T$ can be written $T := \hat{\theta} \circ \bar{S}$, with the $\hat{\theta}$ and $\bar{S}$ defined in condition $M 2$ of Theorem \ref{thm:main}.
These properties, in conjunction with hypotheses $M1-3$ of Theorem \ref{thm:main}, directly imply that condition $(a)$ of Proposition \ref{prop:11} is~verified.

The steps to verify conditions $(b)$ and $(c)$ differ from those in the proof of \cite{fort2003convergence}. Let $\Tilde{S}_n(\theta_n) = \int_z S(z) \Tilde{p}_{\theta, n}(z)\mu(dz)$ be the approximated E step in the Stable Approximate EM, such that $F_n = \hat{\theta} \circ \Tilde{S}_n$. By using uniform continuity properties on compacts, we prove that the following condition is sufficient to verify both $(a)$ and $(b)$:
\begin{equation} \label{eq:sketch_propo_11_sufficient_condition}
 \forall \text{ compact } K, \; \sup{\theta \in K} \, \norm{\Tilde{S}_n(\theta) - \bar{S}(\theta)} \underset{n \to \infty}{\longrightarrow} 0 \, .
\end{equation}

Then, upon replacing $\Tilde{S}_n$ and $\bar{S}$ by their integral forms, it becomes clear that each of the hypotheses \eqref{eq:first_sufficient_condition} or \eqref{eq:second_sufficient_condition} of Theorem \ref{thm:main} is sufficient to verify \eqref{eq:sketch_propo_11_sufficient_condition}. In the end, all conditions are verified to apply Proposition \ref{prop:11}.\\
\\
\textbf{Applying Proposition \ref{prop:11}.} The application of Proposition \ref{prop:11} to the Stable Approximate EM tells us that with probability 1:

$\underset{n \to \infty}{\textnormal{lim sup}}\, j_n < \infty$ and $\brace{\theta_n}_{n\in \N}$ is a compact sequence.

Which is specifically the result $(i)(a)$ of Theorem \ref{thm:main}.

\vspace{12pt}
\noindent \textbf{Verifying the conditions of Proposition \ref{prop:9}.} 
We already have that the likelihood $g$ is a $C^0$, Lyapunov function relative to $(T, \L)$. Thanks to Proposition \ref{prop:11}, we have that $K := Cl\parent{\brace{\theta_n}_n}$ is a compact. Then, $\L \cap K$ is also compact thanks to hypothesis $M3$. Moreover, by definition, the EM sequence verifies: $\brace{\theta_n}_n \in K^{\N}$. The last condition that remains to be shown to apply Proposition \ref{prop:9} is that:
\begin{equation*}
 \underset{n \to \infty}{\lim} | g(\theta_{n+1}) - g \circ T(\theta_n) | = 0 \, .
\end{equation*}

If we apply $(c)$ of Proposition \ref{prop:11} with $K = Cl\parent{\brace{\theta_n}_n}$, we get an almost identical result, but with $\theta_{n+1}$ replaced by $F_n(\theta_n)$. However, $F_n(\theta_n) \neq \theta_{n+1}$ only when $j_{n+1} = j_n + 1$:
\begin{equation*}
 | g(\theta_{n+1}) \!-\! g\!\circ\! T(\theta_n) | \!= | g(\theta_0) \!-\! g\!\circ\! T(\theta_n) | \mathds{1}_{j_{n+1} = j_n + 1} + \! | g \!\circ\! F_n(\theta_n) \!-\! g\!\circ\! T(\theta_n) | \mathds{1}_{j_{n+1} = j_n} \,.
\end{equation*}

We have proven with Proposition \ref{prop:11} that there is only a finite number of such increments. Hence, when $n$ is large enough, $F_n(\theta_n) = \theta_{n+1}$ always, and we have the desired result.\\
\\
\textbf{Applying Proposition \ref{prop:9}}
Since we verify all we need to apply the conclusions of Proposition~\ref{prop:9}:
\begin{itemize}
 \item $\brace{g(\theta_n)}_{n \in \N}$ converges towards a connected component of $g(\L \cap Cl(\brace{\theta_n}_n)) \subset g(\L)$.
 \item If $g(\L \cap Cl(\brace{\theta_n}_n))$ has an empty interior, then the sequence $\brace{g(\theta_n)}_{n \in \N}$ converges towards a $g^* \in \R$ and $\brace{\theta_n}_n$ converges towards the set $\L_{g^*} \cap Cl(\brace{\theta_n}_n) \subseteq \L_{g^*}$.
\end{itemize}

Which are both respectively exactly $(i)(b)$ and $(ii)$ of Theorem \ref{thm:main} and conclude the proof of the Theorem.}

\subsection{Riemann Approximation EM} \label{sect:riemann_em}

\subsubsection{Context and Motivation}
In this section, we introduce one specific case of Approximate EM useful in practice: approximating the conditional probability density function $p_{\theta}(z)$ at the E step by a Riemann sum, in the scenario where the latent variable $z$ is continuous and bounded. We call this procedure the ``Riemann approximation EM''. After motivating this approach, we prove that it is an instance of the Approximate EM algorithm and verifies the hypotheses of Theorem \ref{thm:main}, therefore benefits from the convergence guarantees.

{Consider the case where $z$ is a continuous variable and its} conditional probability $p_{\theta}(z)$ is a continuous function. Even when $h(z; \theta)$ can be computed point by point, a closed form may not exist for the renormalisation term $g(\theta) = \int_z h(z; \theta) dz $. In that case, this integral is usually approximated stochastically with a Monte Carlo estimation, see for instance~\cite{delyon1999convergence, fort2003convergence, allassonniere2021new}. When the dimension is reasonably small, a deterministic approximation through Riemann sums can also be performed. {For the user this can be a welcome simplification, since MCMC methods have a high hyper-parameter complexity (burn-in duration, gain decrease sequence, Gibbs sampler, etc.), whereas the Riemann approximation involves only the position of the Riemann intervals. The choice of which is very guided by the well known theories of integration (Lagrange, Legendre, etc.), and demonstrated in our experiments to have little~impact.}

We introduce the Riemann approximation as a member of the Approximate EM class. Since $z$ is supposed bounded in this section, without loss of generality, we will assume that $z$ is a real variable and $z\in [0,1]$. We recall that $p_{\theta}(z) = h(z; \theta)/g(\theta) = h(z; \theta)/\int_z h(z; \theta) dz $. Instead of using the exact joint likelihood $h(z; \theta)$, we define a sequence of step functions $\brace{\tilde{h}_n}_{n \in \mathbb{N^*}}$ as: $\tilde{h}_n(z; \theta):= h( \lfloor \varphi(n)z \rfloor / \varphi(n); \theta)$. Where $\varphi$ is an increasing function from $\mathbb{N^*} \to \mathbb{N^*}$ such that $\varphi(n) \underset{n \to \infty}{\longrightarrow} \infty$. For the sake of simplicity, we will take $\varphi = Id$, hence $\tilde{h}_n(z; \theta)= h(\lfloor n z \rfloor / n; \theta)$. The following result, however, can be applied to any increasing function $\varphi$ with $\varphi(n) \underset{n \to \infty}{\longrightarrow} \infty$.

With these steps functions, the renormalising factor $\tilde{g}_n(\theta) := \int_z \tilde{h}_n(z; \theta) dz$ is now a finite sum. That is: $\tilde{g}_n(\theta) = \frac{1}{n} \sum_{k=0}^{n-1} h( \lfloor k z \rfloor / n; \theta)$. The approximate conditional probability $\tilde{p}_n(\theta)$ is then naturally defined as: $\tilde{p}_n(\theta) := \tilde{h}_n(z; \theta) / \tilde{g}_n(\theta)$. Thanks to the replacement of the integral by the finite sum, this deterministic approximation is much easier to compute than the real conditional probability.

\subsubsection{Theorem and Proof}
We state and prove the following Theorem for the convergence of the EM with a Riemann approximation. 
\begin{Theorem} \label{thm:riemann_EM}
 {Assume that the hidden variable $z$ is continuous and bounded. Consider the approximation
 \begin{equation*}
 \Tilde{p}_{n, \theta}(z) := \frac{h(\lfloor n z \rfloor / n; \theta)}{\int_{z'} h(\lfloor n z' \rfloor / n; \theta) dz'}\, .
 \end{equation*} 
 We call ``Riemann approximation EM'' the associated Stable Approximate EM. Under conditions $M1-3$ of Theorem \ref{thm:main}, if $z \mapsto S(z)$ is continuous, then the conclusions of Theorem \ref{thm:main} hold for the Riemann approximation EM.}
\end{Theorem}

\begin{proof} {Under the current assumptions, it is sufficient to verify condition \eqref{eq:first_sufficient_condition} in order to apply Theorem \ref{thm:main}. Without loss of generality, we will assume that the bounded variable $z$ is contained in $[0,1]$. Since $S$ is continuous, the first part of the condition is easily verified: $\int_{z=0}^1 S^2_u(z) dz < \infty$.

For the second part of the condition, we consider a compact $K \subseteq \Theta$.} First, note that $h(z; \theta) = \exp(\psi(\theta) + \dotprod{S(z), \phi(\theta)})$ is continuous in $(z, \theta)$, hence uniformly continuous on the compact set $\brack{0,1} \times K$. Additionally, we have: 
\begin{equation*}
\begin{split}
 0 < m := \min{(z, \theta) \in \brack{0,1} \times K} h(z; \theta) \leq h(z; \theta) \, ,\\
 h(z; \theta) \leq \max{(z, \theta) \in \brack{0,1} \times K} h(z; \theta) =: M < \infty \,.
\end{split}
\end{equation*}
where $m$ and $M$ are constants independent of $z$ and $\theta$. This also means that $m \leq g(\theta)=\int_{z=0}^1 h(z;\theta) \leq M$. Moreover, since $\tilde{h}_n(z; \theta) = h\parent{{\floor{n z}}/{n}; \theta}$, then we also have $\forall z \in \brack{0,1}, \theta \in K, n \in \mathbb{N},\; m \leq \tilde{h}_n(z; \theta) \leq M$ and $m\leq \tilde{g}_n(\theta) = \int_{z=0}^1 \tilde{h}_n(z; \theta) \leq M$.

As $h$ is uniformly continuous, $\forall \epsilon > 0, \exists \delta > 0, \forall (z, z') \in \brack{0,1}^2, (\theta, \theta') \in K^2$: 
\begin{equation*}
 \det{z\!-\!z'} \leq \delta \text{ and } \norm{\theta\!-\!\theta'} \leq \delta \!\Rightarrow\! \det{h(z; \theta) \!-\! h(z'; \theta')} \leq \epsilon \, .
\end{equation*}

By definition, ${\floor{n z}}/{n} - z \leq {1}/{n}$. Hence $\exists N \in \mathbb{N}, \forall n \geq N, \floor{n z}/ n - z \leq \delta$. As a consequence: 
\begin{equation*}
 \forall \epsilon > 0,\, \exists N \in \mathbb{N},\, \forall n \geq N,\, \forall (z, \theta) \in \brack{0,1} \times K,\quad \det{h(z; \theta) - \tilde{h}_n(z; \theta)} \leq \epsilon \, .
\end{equation*}

In other words, $\{\tilde{h}_n\}_n$ converges uniformly towards $h$. {Let us fix $\epsilon$ and $N$ now. Then, $\forall n \geq N,\; \forall (z, \theta) \in \brack{0,1}\times K$:}
\begin{equation*}
 \begin{split}
 \Tilde{p}_{\theta, n}(z) - p_{\theta}(z) &= \frac{\tilde{h}_n(z; \theta)}{\int_z \tilde{h}_n(z; \theta) dz} - \frac{h(z; \theta)}{\int_z h(z; \theta) dz} \\
 &=\frac{\tilde{h}_n(z; \theta) - h(z; \theta)}{\int_z \tilde{h}_n(z; \theta) dz} + h(z; \theta) \frac{\int_z \parent{h(z; \theta) - \tilde{h}_n(z; \theta)} dz}{\int_z h(z; \theta)dz \int_z \tilde{h}_n(z; \theta)dz } \\
 &\leq \frac{\epsilon}{m} + M\frac{\epsilon}{m^2}\\
 &= \epsilon \frac{m+M}{m^2} \, .
 \end{split}
\end{equation*}

{
Hence, for this $\epsilon, \; \forall n \geq N:$
\begin{equation*}
 \sup{\theta \in K}\, \int_{z=0}^1 \parent{ \Tilde{p}_{\theta, n}(z) - p_{\theta}(z) }^2 dz \leq \epsilon^2 \parent{\frac{m+M}{m^2}}^2 \, ,
\end{equation*}
which, by definition, means that
\begin{equation*}
 \sup{\theta \in K}\, \int_{z=0}^1 \parent{ \Tilde{p}_{\theta, n}(z) - p_{\theta}(z) }^2 dz \underset{n \to \infty}{\longrightarrow} 0 \, .
\end{equation*}

With this, condition \eqref{eq:first_sufficient_condition} is fully verified, and the conclusions of Theorem \ref{thm:main} are applicable, which concludes the proof.}
\end{proof}

\subsection{Tempered EM} \label{sect:tmp_em}

\subsubsection{Context and Motivation}
{In this section, we consider another particular case of Deterministic Approximate EM: the Tempered EM (or ``tmp-EM''), first introduced in \cite{ueda1998deterministic}. We first prove that under mild conditions, tmp-EM verifies the hypotheses of Theorem \ref{thm:main}, hence benefits from the state of the art EM convergence guarantees. In particular, we prove that the choice of the temperature profile is almost completely free. This justifies the use of a wider array of temperature profiles than the ones specified in \cite{ueda1998deterministic,naim2012convergence}. Then, we demonstrate the interest of tmp-EM with non-monotonous profiles, on Mixture estimation problems with hard to escape initialisations.

\textit{Tempering} or \textit{annealing} is a common technique in the world of non-convex optimisation. With a non-convex objective function $g$, naively following the gradients usually leads to undesirable local extrema close to the initial point. A common remedy is to elevate $g$ to the power $\frac{1}{T_n}$, with $\{T_n\}_{n\in \mathbb{N}}$ a sequence of temperatures tending towards 1 as the number $n$ of steps of the procedure increases. When $T_n > 1$, the shape of the new objective function $g^{\frac{1}{T_n}}$ is flattened, making the gradients less strong, and the potential wells less attractive, but without changing the hierarchy of the values. This allows the optimisation procedure to explore more before converging. This concept is put in application in many state of the art procedures. The most iconic probably being the Simulated Annealing, introduced and developed in \cite{kirkpatrick1983optimization, van1987simulated, aarts1988simulated}, where in particular $T_n \longrightarrow 0$ instead of 1. It is one of the few optimisation technique proven to find global optimum of non-convex functions. The Parallel Tempering (or Annealing MCMC) developed in \cite{swendsen1986replica, geyer1995annealing, hukushima1996exchange} also makes use of these ideas to improve the MCMC simulation of a target probability distribution.

Since non-convex objectives are a common occurrence in the EM literature, \cite{ueda1998deterministic} introduced the \textit{Deterministic Annealed EM}, where, in the E step, the conditional probability is replaced by a tempered approximated distribution: $\Tilde{p}_{n, \theta}(z) \propto p_{\theta}^{\frac{1}{T_n}}(z) \propto h(z; \theta)^{\frac{1}{T_n}}$ (renormalised such that $\int_z \Tilde{p}_{n, \theta}(z) dz = 1$). For the temperature profile they consider specific sequences $\{T_n\}_{n\in \mathbb{N}}$ decreasingly converging to 1. Another specific, non-monotonous, temperature scheme was later proposed by \citet{naim2012convergence}. In both cases, theoretical convergence guarantees are lacking. Later, in \cite{allassonniere2021new}, tempering has been applied to the SAEM, with convergence guarantees provided for any temperature scheme $T_n \longrightarrow 1$.

With Theorem \ref{thm:tempered_EM}, we show that the Deterministic Anealing EM, or tmp-EM, is a specific case of the Deterministic Approximate EM of Section \ref{sect:approx_EM}, hence can benefit from the same convergence properties. In particular, we show that any temperature scheme $\{T_n\}_n \in (\R^*)^{\mathbb{N}}, \: T_n \underset{n \to \infty}{\longrightarrow} 1$ guarantees the state of the art convergence.}
\begin{Remark}
 Elevating $p_{\theta}(z)$ to the power $\frac{1}{T_n}$, as is done here and in \cite{ueda1998deterministic, naim2012convergence}, is not equivalent to elevating to the power $\frac{1}{T_n}$ the objective function $g(\theta)$, which would be expected for a typical annealed or tempered optimisation procedure. It is not equivalent either to elevating to the power $\frac{1}{T_n}$ the proxy function $\mathbb{E}_{z \sim p_{\theta_n}(z)} \brack{ h(z; \theta) }$ that is optimised in the M step. Instead, the weights $p_{\theta_n}(z)$ (or equivalently, the terms $h(z; \theta_n)$) used in the calculation of $\mathbb{E}_{z \sim p_{\theta_n}(z)} \brack{ h(z; \theta) }$ are the tempered terms. This still results in the desired behaviour and is only a more ``structured'' tempering. Indeed, with this tempering, it is the estimated distribution of the latent variable $z$ that are made less unequivocal, with weaker modes, at each step. This forces the procedure to spend more time considering different configurations for those variables, which renders as a result the optimised function $\mathbb{E}_{z \sim p_{\theta_n}(z)} \brack{ h(z; \theta) }$ more ambiguous regarding which $\theta$ is the best, just as intended. Then, when $n$ large, the algorithm is allowed to converge, as $T_n \longrightarrow 1$ and $\mathbb{E}_{z \sim \tilde{p}_{n, \theta}(z)} \longrightarrow \mathbb{E}_{z \sim p_{\theta}(z)}$.
\end{Remark}

\subsubsection{Theorem}
We now provide the convergence Theorem for the Approximate EM with the tempering approximation. In particular, this result highlights that there are almost no constraints on the temperature profile to achieve convergence. 
\begin{Theorem} \label{thm:tempered_EM}
 {Let $T_n$ be a sequence of non-zero real numbers. Consider the approximation introduced in \cite{ueda1998deterministic}:
 \begin{equation*}
 \Tilde{p}_{n, \theta}(z) := \frac{p_{\theta}^{\frac{1}{T_n}}(z)}{\int_{z'} p_{\theta}^{\frac{1}{T_n}}(z') dz'} \, .
 \end{equation*}
 We call ``Tempered EM'' the associated Stable Approximate EM. Define $\overline{\mathcal{B}}(1, \epsilon)$ be the closed ball centred in 1 and with radius $\epsilon \in \R_+$. Under conditions $M1-3$ of Theorem \ref{thm:main}, if $T_n \underset{n \to \infty}{\longrightarrow} 1$ and for any compact $K \subseteq \Theta, \; \, \exists \epsilon \in ]0,1[, \; \, \forall \alpha \in \overline{\mathcal{B}}(1, \epsilon)$:
 \begin{itemize}
 \item[] (T1) $\ \sup{\theta \in K}\, \int_z p_{\theta}^{\alpha}(z) dz < \infty \,$, 
 \item[] (T2) $\ \forall u \in \llbracket{1, q\rrbracket}, \quad \sup{\theta \in K}\, \int_z S^2_u(z) p_{\theta}^{\alpha}(z) dz < \infty \,$,
 \end{itemize}
 then the conclusions of Theorem \ref{thm:main} hold for the Tempered EM.}
\end{Theorem}

\begin{Remark}
 {The added condition that $(T1)$ and $(T2)$ must hold for all $\alpha$ in a ball $\overline{\mathcal{B}}(1, \epsilon)$ is very mild. Indeed, in Section \ref{sect:theoretical_examples}, we show classical examples that easily verify the much stronger condition that $(T1)$ and $(T2)$ hold for all $\alpha \in \R^*_+$.}
\end{Remark}

\subsubsection{Sketch of Proof} 
{The detailed proof of Theorem \ref{thm:tempered_EM} can be found in Appendix \ref{sect:proofs_tempered}, we propose here an abbreviated version where we highlight the key steps.\\
Under the current assumptions, it is sufficient to verify the second part of condition \eqref{eq:second_sufficient_condition} in order to apply Theorem \ref{thm:main}. To that end, we must control the integral
\begin{equation*}
 \int_z \parent{ \frac{\Tilde{p}_{\theta, n}(z)}{p_{\theta}(z)} - 1 }^2 p_{\theta}(z) dz \, ,
\end{equation*}
for all $\theta$ in a compact $K \subseteq \Theta$. First, with a Taylor development in the term $\parent{\frac{1}{T_n}-1}$, which converges toward 0 when $n \to \infty$, we control the difference $\parent{\Tilde{p}_{\theta, n}(z) - p_{\theta}(z)}^2$:}
\begin{equation*}
 \Bigg(\!\frac{p_{\theta}(z)^{\frac{1}{T_n}}}{\int_{z'} p_{\theta}(z')^{\frac{1}{T_n}}}\! -\! p_{\theta}(z) \!\Bigg)^2 \leq 2\parent{\frac{1}{T_n}-1}^2 p_{\theta}(z)^2 \!\Bigg(\!\!\Big(\ln p_{\theta}(z)\, e^{a(z, \theta, T_n)}\Big)^2 \!A(\theta, T_n) \!+\! B(\theta, T_n) \!\! \Bigg) \, .
\end{equation*}
where the terms $A(\theta, T_n)$, $B(\theta, T_n)$ and $a(z, \theta, T_n)$ come from the Taylor development. Then, we can control the integral of interest:
\begin{equation} \label{eq:thm2_sketch_proof_inequality}
\begin{split}
 \int_z \frac{\parent{ \frac{p_{\theta}(z)^{\frac{1}{T_n}}}{\int_{z'} p_{\theta}(z')^{\frac{1}{T_n}}} - p_{\theta}(z) }^2}{p_{\theta}(z)} dz \leq & 2 \parent{\!\frac{1}{T_n}\!-\!1}^2\! A(\theta, T_n)\! \int_z\! p_{\theta}(z) e^{2a(z, \theta, T_n)} \ln^2 p_{\theta}(z)dz \\
 &+ 2 \parent{\frac{1}{T_n}-1}^2 B(\theta, T_n) \, .
\end{split}
\end{equation}

From the properties of the Taylor development, we prove that $A(\theta, T_n)$ and $B(\theta, T_n)$ both have upper bounds involving only $\int_z p_{\theta}(z) \ln p_{\theta}(z)$ and $\int_z p_{\theta}(z)^{\frac{1}{T_n}} \ln p_{\theta}(z)$. In a similar fashion, the term $\int_z p_{\theta}(z) e^{2a(z, \theta, T_n)} \ln^2 p_{\theta}(z)$ has an upper bound involving only $\int_z p_{\theta}(z) \ln^2 p_{\theta}(z) dz$ and $\int_z p_{\theta}(z)^{\frac{2}{T_n}-1} \ln^2 p_{\theta}(z) dz$.

Using the hypotheses of the Theorem, we prove that for any $\alpha \in \overline{\mathcal{B}}(1, \epsilon)$ and $\theta \in K$ the terms $\int_z p_{\theta}(z)^{\alpha} \ln p_{\theta}(z) $ and $\int_z p_{\theta}(z)^{\alpha} \ln^2 p_{\theta}(z) $ are both upper bounded by a constant $C$ independent of $\theta$ and $\alpha$.

Since $T_n \underset{n \to \infty}{\longrightarrow} 1$, then when $n$ is large enough, $\frac{1}{T_n} \in \overline{\mathcal{B}}(1, \epsilon)$ and $\frac{2}{T_n}-1 \in \overline{\mathcal{B}}(1, \epsilon)$. Hence, the previous result applies to the upper bounds of $A(\theta, T_n)$, $B(\theta, T_n)$ and $\int_z p_{\theta}(z) e^{2a(z, \theta, T_n)}$ $ \ln^2 p_{\theta}(z)dz$. As a result, these three terms are respectively upper bounded by $C_1$, $C_2$ and $C_3$, three constants independent of $\theta$ and $T_n$.

The inequality \eqref{eq:thm2_sketch_proof_inequality} then becomes:
\begin{equation*}
 \int_z \frac{1}{p_{\theta}(z)} \parent{ \frac{p_{\theta}(z)^{\frac{1}{T_n}}}{\int_{z'} p_{\theta}(z')^{\frac{1}{T_n}}} - p_{\theta}(z) }^2 dz \leq 2 \parent{\frac{1}{T_n}-1}^2 C_1 C_2 + 2 \parent{\frac{1}{T_n}-1}^2 C_3 \, .
\end{equation*}

By taking the supremum in $\theta \in K$ and the limit when $n \longrightarrow \infty$, we get the desired result:
\begin{equation*}
 \sup{\theta \in K} \int_z \frac{1}{p_{\theta}(z)} \parent{ \frac{p_{\theta}(z)^{\frac{1}{T_n}}}{\int_{z'} p_{\theta}(z')^{\frac{1}{T_n}}} - p_{\theta}(z) }^2 dz \underset{n \to \infty}{\longrightarrow} 0 \, .
\end{equation*}

With condition \eqref{eq:second_sufficient_condition} verified, the conclusions of Theorem \ref{thm:main} are applicable, which concludes the proof.


\subsubsection{Examples of Models That Verify the Conditions} \label{sect:theoretical_examples}
In this section, we illustrate that the conditions of Theorem \ref{thm:tempered_EM} are easily met by common models. We take two examples, first the Gaussian Mixture Model (GMM) where the latent variables belong to a finite space, then the Poisson count with random effect, where the latent variables live in a continuous space. {As mentioned, these examples are shown to not only verify all hypotheses of Theorem \ref{thm:tempered_EM}, but also to verify $(T1)$ and $(T2)$ for any $\alpha \in \R_+^*$.}

\paragraph{Gaussian Mixture Model}
Despite being one of the most common models the EM is applied to, the GMM have many known irregularities and pathological behaviours, see \cite{titterington_85, ho2016convergence}. Although some recent works, such as \cite{dwivedi2020singularity, dwivedi2020sharp}, tackled the theoretical analysis of EM for GMM, none of the convergence results for the traditional EM and its variants proposed by \cite{wu1983convergence, lange1995gradient, delyon1999convergence, fort2003convergence} apply to the GMM. The hypothesis that the GMM fail to verify is the condition that the level lines have to be compact ($M2$ $(d)$ in this paper). All is not lost however for the GMM, indeed, the model verifies all the other hypotheses of the general Theorem \ref{thm:main} as well as the tempering hypotheses of Theorem \ref{thm:tempered_EM}. Moreover, in this paper as in the others, the only function of the unverified hypothesis $M2$ $(d)$ is to ensure in the proof that the EM sequence stays within a compact. The latter condition is the actual relevant property to guarantee convergence at this stage of the proof. This means that, in practice, if an tempered EM sequence applied to a GMM is observed to remain within a compact, then all the conditions for convergence are met, Theorem \ref{thm:tempered_EM} applies, and the sequence is guaranteed to converge towards a critical point of the likelihood function.

In the following, we provide more details on the GMM likelihoods and the theorem hypotheses they verify. First, note that the GMM belongs to the exponential family with the complete likelihood:
\begin{equation}\label{eq:MoG_complete_likelihood}
\begin{split}
 h(z; \theta) = \prod_{i=1}^N \sum_{k=1}^K exp\Bigg( \frac{1}{2}\Big( &-(x^{(i)}-\mu_k)^T \Theta_k (x^{(i)} -\mu_k)\\
 &+ \ln(\det{\Theta_k}) + 2 \ln (\pi_k) - p \ln (2\pi) \Big) \Bigg){\mathds{1}_{z^{(i)}=k}} \, ,
\end{split}
\end{equation}
and the observed likelihood:
\begin{equation}\label{eq:MoG_likelihood}
\begin{split}
 g(\theta) = \prod_{i=1}^N \sum_{k=1}^K exp\Bigg( \frac{1}{2}\Big( &-(x^{(i)}-\mu_k)^T \Theta_k (x^{(i)} -\mu_k)\\
 &+ \ln (\det{\Theta_k})+ 2 \ln (\pi_k) - p \ln (2\pi) \Big) \Bigg)\, .
\end{split}
\end{equation}

This is an exponential model with parameter:
\begin{equation*}
 \theta := \parent{ \brace{\pi_k}_{k=1}^K, \brace{\mu_k}_{k=1}^K, \brace{\Theta_k}_{k=1}^K } \in \Big\{ \brace{\pi_k}_{k} \in \brack{0, 1}^K \Bigg| \sum_k \pi_k = 1 \Big\} \otimes \R^{p \times K} \otimes {S_p^{++}}^{K} . \\
\end{equation*}
where $S_p^{++}$ is the cone of symmetric positive definite matrices of size $p$. The verification of conditions \emph{M1--3} for the GMM (except $M2$ $(d)$ of course) is a classical exercise since these are the conditions our Theorem shares with any other EM convergence result on the exponential family. We focus here on the hypotheses specific to our Deterministic Approximate EM.\\
\\
\textit{Condition on \texorpdfstring{$\int_z p_{\theta}^{\alpha}(z) dz $}{TEXT}}. Let $\alpha \in \R^*_+$, in the finite mixture case, the integrals on $z$ are finite sums: 
\begin{equation*}
 \int_z p_{\theta}^{\alpha}(z) dz = \sum_{k} p_{\theta}^{\alpha}(z=k) \, .
\end{equation*}

Which is continuous in $\theta$ since $\theta \mapsto p_{\theta}(z=k) = h(z=k; \theta) / g(\theta)$ is continuous. Hence
\begin{equation*}
 \forall \alpha \in \R^*_+, \quad \sup{\theta \in K}\, \int_z p_{\theta}^{\alpha}(z) dz < \infty \, .
\end{equation*}
\noindent
\textit{Condition on \texorpdfstring{$\int_z S^2_u(z) p_{\theta}^{\alpha}(z) dz$}{TEXT}}. The previous continuity argument is still valid.

\paragraph{Poisson Count with Random Effect}
This model is discussed in \cite{fort2003convergence}, the authors prove, among other things, that this model verifies the conditions $M1-3$. Here is a brief description of the model: the hidden variables $\brace{z^{(i)}}_{i=1}^N$ are distributed according to an autoregressive process of order 1: $z^{(i)} := a z^{(i-1)} + \sigma \epsilon_i$, with $|a|<1$, $\sigma>0$ and the $\brace{\epsilon^{(i)}}_{i=1}^N$ are iid standard gaussian. Conditionally to $\brace{z^{(i)}}_{i=1}^N$, the observed variables $x^{(i)}$ are independent Poisson variables with parameter $\exp(\theta+z^{(i)})$. The complete likelihood of the model, not accounting for irrelevant constants,~is:
\begin{equation}\label{eq:Poisson_complete_likelihood}
 h(z; \theta) = e^{\theta \sum_{i=1}^N x^{(i)}}. exp\parent{-e^{\theta} \sum_{i=1}^N e^{z^{(i)}}} \, .
\end{equation}
$g(\theta) = \int_z h(z; \theta) dz $ can be computed analytically up to a constant:
\begin{equation}\label{eq:Poisson_likelihood}
\begin{split}
 g(\theta) &= \int_{z \in \R^N} h(z; \theta) dz \\
 &= e^{\theta \sum_{i=1}^N x^{(i)}} \int_{z \in \R^N} exp\parent{-e^{\theta} \sum_{i=1}^N e^{z^{(i)}}} dz \\
 &= e^{\theta \sum_{i=1}^N x^{(i)}} \prod_{i=1}^N \int_{z^{(i)} \in \R} exp\parent{-e^{\theta} e^{z^{(i)}}} dz^{(i)} \\
 &= e^{\theta \sum_{i=1}^N x^{(i)}} \parent{\int_{u \in \R_+} \frac{exp\parent{-u}}{u} du }^N\\
 &= e^{\theta \sum_{i=1}^N x^{(i)}} E_1(0)^N \, ,
\end{split}
\end{equation}
where $E_1(0)$ is a finite, non zero, constant, called ``exponential integral'', in particular independent of $\alpha$ and $\theta$.\\
\\
\textit{Condition on \texorpdfstring{$\int_z p_{\theta}^{\alpha}(z) dz $}{TEXT}}. Let $K$ be a compact in $\Theta$.\\
We have $p_{\theta}(z) = \frac{h(z; \theta)}{g(\theta)}$. Let us compute $\int_z h(z; \theta)^{\alpha} $ for any positive $\alpha$. The calculations work as in Equation~\eqref{eq:Poisson_likelihood}:
\begin{equation*}
 \begin{split}
 \int_{z \in \R^N} h(z; \theta)^{\alpha} &= e^{\alpha \theta \sum_{i=1}^N x^{(i)}} \prod_{i=1}^N \int_{z^{(i)} \in \R} exp\parent{-\alpha e^{\theta} e^{z^{(i)}}} dz^{(i)} \\
 &= e^{\alpha \theta \sum_{i=1}^N x^{(i)}} E_1(0)^N \, .
 \end{split}
\end{equation*}
Hence:
\begin{equation*}
 \int_z p_{\theta}^{\alpha}(z) dz = E_1(0)^{(1-\alpha) N}\, .
\end{equation*}

Since $E_1(0)$ is finite, non zero, and independent of $\theta$, we easily have:
\begin{equation*}
 \forall \alpha \in \R^*_+, \quad \sup{\theta \in K}\, \int_z p_{\theta}^{\alpha}(z) dz < \infty \, .
\end{equation*}
$\theta$ does not even have to be restricted to a compact.\\
\\
\textit{Condition on \texorpdfstring{$\int_z S^2_u(z) p_{\theta}^{\alpha}(z) dz$}{TEXT}}. Let $K$ be a compact in $\Theta$ and $\alpha$ a positive real number.\\
In this Poisson count model, $S(z) = \sum_{i=1}^N e^{z^{(i)}} \in \R$. We have:
\begin{equation} \label{eq:Poisson_integral_S2}
 S^2(z) p_{\theta}^{\alpha}(z) = \parent{\sum_{i=1}^N e^{z^{(i)}}}^2 \frac{exp\parent{-\alpha e^{\theta} \sum_{i=1}^N e^{z^{(i)}}}}{E_1(0)^{\alpha N}} \, .
\end{equation}

First, let us prove that the integral is finite for any $\theta$. We introduce the variables $u_k := \sum_{l=1}^k e^{z_l}$. The Jacobi matrix is triangular and its determinant is $ \prod_k e^{z^{(i)}} = \prod_k u_k$.
\begin{equation*}
\begin{split}
 \int_z\! S^2(z) p_{\theta}^{\alpha}(z) dz \!&= \!\frac{\int_{z\in \R^N}\! \parent{\sum_k \!e^{z^{(i)}}}^2\! exp\parent{\!-\alpha e^{\theta} \!\sum_k\! e^{z^{(i)}}\!} \!dz}{E_1(0)^{\alpha N}}\, .
\end{split}
\end{equation*}

Which is proportional to:
\begin{equation*}
 \int_{u_1 = 0}^{+\infty} \!u_1 \int_{u_2 = u_1}^{+\infty} \!u_2... \int_{u_N = u_{N-1}}^{+\infty} \!u_N^3 \, e^{-\alpha e^{\theta} u_N}\, du_N...du_2 du_1 \, .
\end{equation*}
where we removed the finite constant $\frac{1}{E_1(0)^{\alpha N}}$ for clarity.
This integral is finite for any $\theta$ because the exponential is the dominant term around $+\infty$. Let us now prove that $\theta \mapsto \int_z S^2(z) p_{\theta}^{\alpha}(z) dz$ is continuous. From Equation~\eqref{eq:Poisson_integral_S2}, we have that
\begin{itemize}
 \item $z \mapsto S^2(z) p_{\theta}^{\alpha}(z)$ is measurable on $\R^N$.
 \item $\theta \mapsto S^2(z) p_{\theta}^{\alpha}(z)$ is continuous on $K$ (and on $\Theta = \R$). 
 \item With $\theta_M := \underset{\theta \in K}{min} \, \theta$, then $\forall \theta \in K, \; 0 \leq S^2(z) p_{\theta}^{\alpha}(z) \leq S^2(z) p_{\theta_M}^{\alpha}(z) $ 
\end{itemize}

Since we have proven that $S^2(z) p_{\theta_M}^{\alpha}(z) < \infty$, then we can apply the intervertion Theorem and state that $\theta \mapsto \int_z S^2(z) p_{\theta}^{\alpha}(z) dz$ is continuous.\\
It directly follows that:
\begin{equation*}
 \forall \alpha \in \R^*_+, \quad \sup{\theta \in K}\, \int_z S^2(z) p_{\theta}^{\alpha}(z) dz < \infty \, .
\end{equation*}

Note that after the change of variable, the integral could be computed explicitly, but involves $N$ successive integration of polynomial $\times$ exponential function products of the form $P(x) e^{-\alpha e^{\theta} x}$. This would get tedious, especially since after each successful integration, the product with the next integration variable $u_{k-1}$ increases by one the degree of the polynomial, i.e. starting from 3, the degree ends up being $N+2$. We chose a faster path.

\subsection{Tempered Riemann Approximation EM} \label{sect:tmp_riemann_em}

\subsubsection*{Context, Theorem and Proof}
The Riemann approximation of Section \ref{sect:riemann_em} makes the EM computations possible in hard cases, when the conditional distribution has no analytical form for instance. It is an alternative to the many stochastic approximation methods (SAEM, MCMC-SAEM, etc.) that are commonly used in those cases. The tempering approximation of Section \ref{sect:tmp_em} is used to escape the initialisation by allowing the procedure to explore more the likelihood profile before committing to convergence. We showed that both these approximation are particular cases of the wider class of Deterministic Approximate EM, introduced in Section~\ref{sect:approx_EM}. However, since they fulfil different purposes, it is natural to use them in coordination and not as alternatives of one another. In this section, we introduce another instance of the Approximate EM: a combination of the tempered and Riemann sum approximations. This ``tempered Riemann approximation EM'' (tmp-Riemann approximation) can compute EM steps when there is no closed form thanks to the Riemann sums as well as escape the initialisation thanks to the tempering. For a bounded latent variable $z \in \brack{0,1}$, we define the approximation as: $\Tilde{p}_{n, \theta}(z) := h(\lfloor n z \rfloor / n; \theta)^{\frac{1}{T_n}} / \int_{z'} h(\lfloor n z' \rfloor / n; \theta)^{\frac{1}{T_n}} dz'$, for a sequence $\{T_n\}_n \in (\R^*_+)^{\mathbb{N}}, \: T_n \underset{n \to \infty}{\longrightarrow} 1$.

In the following Theorem, we prove that the tempered Riemann approximation EM verifies the applicability conditions of Theorem \ref{thm:main} with no additional hypothesis from the regular Riemann approximation EM covered by Theorem \ref{thm:riemann_EM}. 

\begin{Theorem} \label{thm:tmp_riemann_em}
 Under conditions $M1-3$ of Theorem \ref{thm:main}, and when $z$ is bounded, the (Stable) Approximate EM with $\Tilde{p}_{n, \theta}(z) := \frac{h(\lfloor n z \rfloor / n; \theta)^{\frac{1}{T_n}}}{\int_{z'} h(\lfloor n z' \rfloor / n; \theta)^{\frac{1}{T_n}} dz'}$, which we call ``tempered Riemann approximation EM'', verifies the remaining conditions of applicability of Theorem \ref{thm:main} as long as $z \mapsto S(z)$ is continuous and $\{T_n\}_n \in (\R^*_+) ^{\mathbb{N}}, \: T_n \underset{n \to \infty}{\longrightarrow} 1$.
\end{Theorem}

\begin{proof} 
 This proof of Theorem \ref{thm:tmp_riemann_em} is very similar to the proof of Theorem \ref{thm:riemann_EM} for the regular Riemann approximation EM. The first common element is that for the tempered Riemann approximation EM, the only remaining applicability condition of the general Theorem \ref{thm:main} to prove is also: 
 \begin{equation*}
 \forall \text{compact} \; K \subseteq \Theta, \; \sup{\theta \in K} \, \int_z \parent{ \Tilde{p}_{\theta, n}(z) - p_{\theta}(z) }^2 dz \underset{n \to \infty}{\longrightarrow} 0 \, .
 \end{equation*}
 
 In the proof of Theorem \ref{thm:riemann_EM}, we proved that having the uniform convergence of the approximated complete likelihood $\{\tilde{h}_n\}_n$ towards the real $h$ - with both $\tilde{h}_n(z; \theta)$ and $h(z; \theta)$ uniformly bounded - was sufficient to fulfil this condition. Hence, we prove In this section, that these sufficient properties still hold, even with the tempered Riemann approximation, where $\tilde{h}_n(z; \theta)\! :=\! h\!\parent{{\floor{n z}}/{n}; \theta}^{\frac{1}{T_n}}$.
 
 We recall that $h(z; \theta)$ is uniformly continuous on the compact set $\brack{0,1} \times K$, and verifies: 
 \begin{equation*}
 0 < m \leq h(z; \theta) \leq M < \infty \,.
 \end{equation*}
 where $m$ and $M$ are constants independent of $z$ and $\theta$.
 
 Since $T_n > 0, \; T_n \underset{n \to \infty}{\longrightarrow} 1$, then the sequence $\{1/T_n\}_{n}$ is bounded. Since $\tilde{h}_n(z; \theta) = h\parent{{\floor{n z}}/{n}; \theta}^{\frac{1}{T_n}}$, with $0 < m \leq h\parent{{\floor{n z}}/{n}; \theta} \leq M < \infty$ for any $z, \theta$ and $n$, then we also have:
 \begin{equation*}
 0 < m' \leq \tilde{h}_n(z; \theta) \leq M' < \infty \,,
 \end{equation*}
 with $m'$ and $M'$ constants independent of $z, \theta$ and $n$.
 
 We have seen in the proof of Theorem \ref{thm:riemann_EM}, that:
 \begin{equation*}
 \forall \epsilon > 0,\, \exists N \in \mathbb{N},\, \forall n \geq N,\, \forall (z, \theta) \in \brack{0,1} \times K,\quad \det{h(z; \theta) - h\parent{{\floor{n z}}/{n}; \theta}} \leq \epsilon \, .
 \end{equation*}
 
 To complete the proof, we control in a similar way the difference $h\parent{{\floor{n z}}/{n}; \theta} - h\parent{{\floor{n z}}/{n}; \theta}^{\frac{1}{T_n}}$. The function $(h, T) \in [m, M] \times [T_{min}, T_{max}] \mapsto h^{\frac{1}{T}} \in \R$ is continuous on a compact, hence uniformly continuous in $(h, T)$. As a consequence: $\forall \epsilon > 0, \exists \delta > 0, \forall (h, h') \in \brack{m, M}^2, (T, T') \in [T_{min}, T_{max}]^2,$
 \begin{equation*}
 \det{h\!-\!h'} \leq \delta \text{ and } \det{T\!-\!T'} \leq \delta \implies \det{h^{\frac{1}{T}} \!-\!(h')^{\frac{1}{T'}})} \leq \epsilon \, .
 \end{equation*}
 
 Hence, with $N \in \N$ such that $\forall n \geq N, \det{T_n-1} \leq \delta$, we have:
 \begin{equation*}
 \begin{split}
 \forall n \geq N, &\forall (z, \theta) \in\brack{0,1} \times K,\\
 &\det{h\parent{{\floor{n z}}/{n}; \theta}-h\parent{{\floor{n z}}/{n}; \theta}^{\frac{1}{T_n}}} \leq \epsilon \, .
 \end{split}
 \end{equation*}
 
 In the end, $\forall \epsilon > 0, \exists N \in \mathbb{N}, \forall n \geq N, \forall (z, \theta) \in \brack{0,1} \times K$:
 \begin{equation*}
 \begin{split}
 \det{h(z; \theta) - \tilde{h}_n\parent{z; \theta}} &= \det{h(z; \theta) - h\parent{{\floor{n z}}/{n}; \theta}^{\frac{1}{T_n}}}\\
 &\leq \det{h(z; \theta) - h\parent{{\floor{n z}}/{n}; \theta}} + \det{h\parent{{\floor{n z}}/{n}; \theta}-h\parent{{\floor{n z}}/{n}; \theta}^{\frac{1}{T_n}}}\\
 &\leq 2\epsilon \, .
 \end{split}
 \end{equation*}
 
 In other words, we have the uniform convergence of $\{\tilde{h}_n\}$ towards $h$. From there, we conclude following the same steps as in the proof of Theorem \ref{thm:riemann_EM}. 
\end{proof}

\section{Results}\label{sect:exp}
In this section, we describe experiments that explore each of the three studied methods: Riemann EM, tempered EM and tempered Riemann.

\subsection{Riemann Approximation EM: Two Applications} \label{sect:exp_riemann_em}
\subsubsection{Application to a Gaussian Model with the Beta Prior} \label{sect:riemann_em_example}
We demonstrate the interest of the method on a example with a continuous bounded random variable following a Beta distribution $z \sim Beta(\alpha, 1)$, and an observed random variable following $x \sim \mathcal{N}(\lambda z, \sigma^2)$. In other words, with $\epsilon \sim \mathcal{N}(0,1)$ independent of $z$: 
\begin{equation*}
 x = \lambda z + \sigma \epsilon \, .
\end{equation*}

This results in a likelihood belonging to the exponential family:
\begin{equation*}
 h(z; \theta) =\frac{\alpha z^{\alpha-1}}{\sqrt{2 \pi \sigma^2}} exp\parent{-\frac{(x-\lambda z)^2}{2 \sigma^2}}\, .
\end{equation*}

Since $z$ is bounded, and everything is continuous in the parameter $(\alpha, \lambda, \sigma^2)$, this model easily verifies each of the conditions \emph{M1--3}.
The E step with this model involves the integral $\int_z z^{\alpha} exp\parent{-\frac{(x-\lambda z)^2}{2 \sigma^2}} dz$, a fractional moment of the Gaussian distribution. Theoretical formulas exists for these moments, see \cite{winkelbauer2012moments}, however they involve Kummer’s confluent hypergeometric functions, which are infinite series. Instead, we use the Riemann approximation to run the EM algorithm with this model: $\tilde{h}_n(z; \theta):= h( \lfloor \varphi(n)z \rfloor / \varphi(n); \theta)$. As done previously, we take, without loss of generality, $\varphi(n):=n$ for the sake of simplicity. The E step only involves the $n$ different values taken by the step function probabilities $h( \lfloor n z \rfloor / n; \theta)$:
\begin{equation*}
 \tilde{p}_{\theta, n}^{(i)}\parent{\frac{k}{n}} = \frac{h^{(i)}(\frac{k}{n}; \theta)}{\frac{1}{n} \sum_{l=0}^{n-1} h^{(i)}(\frac{l}{n}; \theta)} \, .
\end{equation*}
where the exponent $(i)$ indicates the index of the observation $x^{(i)}$. To express the corresponding M step in a digest way, let us define the operator $\Psi^{(i)} : \R^{[0, 1]} \longrightarrow \R$ such that, for any $f : [0, 1] \longrightarrow \R$:
\begin{equation*}
 \Psi^{(i)} \circ f = \sum_{k=0}^{n-1} \tilde{p}_{\theta, n}^{(i)} \parent{\frac{k}{n}} \int_{z=k/n}^{(k+1)/n} f(z) dz \, .
\end{equation*} 

Then, the M step can be expressed as:
\begin{equation} \label{eq:gaussian_beta_M_step}
 \begin{split}
 \frac{1}{\hat{\alpha}} &= -\frac{1}{N} \sum_{i=1}^N \Psi^{(i)} \circ \ln(z) \, , \\
 \hat{\lambda} &= \frac{\sum_{i=1}^N \Psi^{(i)} \circ (x^{(i)} z)}{\sum_{i=1}^N \Psi^{(i)} \circ z^2} \, , \\
 \hat{\sigma}^2 &= \frac{1}{N} \sum_{i=1}^N \Psi^{(i)} \circ \parent{x^{(i)}-\hat{\lambda}z}^2 \, .
 \end{split}
\end{equation}
where we took the liberty of replacing $f$ by $f(z)$ in these equations for the sake of simplicity. Here $N$ is the total number of observations: $x := (x^{(1)}, ..., x^{(N)})$ iid.

\begin{figure}[tbhp]
 \subfloat[]{\includegraphics[width=0.5\linewidth, valign=c]{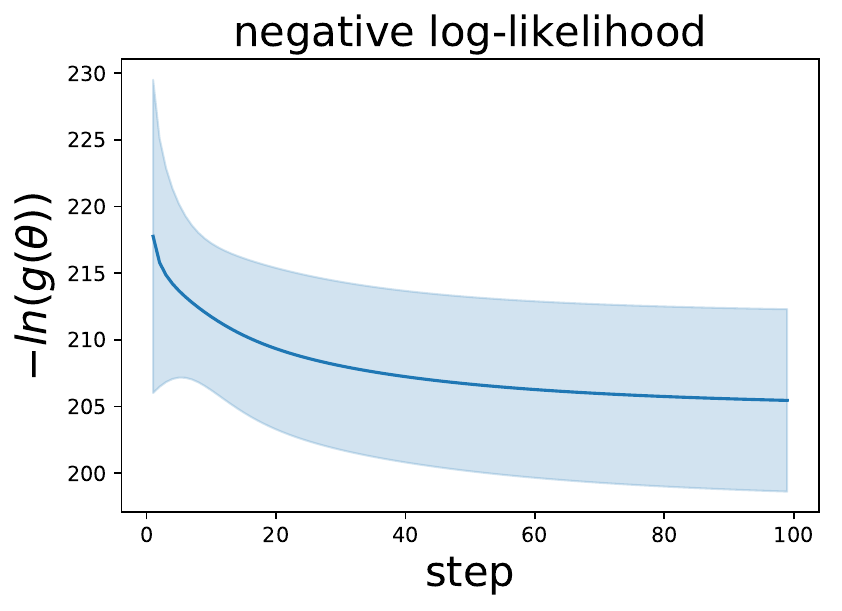}}
 \subfloat[]{
 \adjustbox{valign=c}{\begin{tabular}{cc}
 \toprule
 Metric & Average (std)\\
 \midrule \\
 $\frac{(\alpha - \hat{\alpha})^2}{\alpha^2}$ & 0.129 (0.078) \\
 $\frac{(\lambda - \hat{\lambda})^2}{\lambda^2}$ & 0.048 (0.058)\\
 $\frac{(\sigma - \hat{\sigma})^2}{\sigma^2}$ & 0.108 (0.119) \\ \bottomrule
 \end{tabular}}
 }
\caption{(\textbf{a}). Average values, with standard deviation, over 2000 simulations of the negative log-likelihood along the steps of the Riemann EM. The Riemann EM increases the likelihood. (\textbf{b}). Average and standard deviation of the relative parameter reconstruction errors at the end of the Riemann EM.}
\label{fig:riemann}
\end{figure}

We test this algorithm on synthetic data. With real values $\alpha = 2, \lambda = 5, \sigma = 1.5$, we generate a dataset with $N=100$ observations and run the Riemann EM with random initialisation. This simulation is ran 2000 times. We observe that the Riemann EM is indeed able to increase the likelihood, despite the EM being originally intractable. On \figurename~\ref{fig:riemann}, we display the average trajectory, with standard deviation, of the negative log-likelihood $-\ln(g(\theta))$ during the Riemann EM procedure. The profile is indeed decreasing. The standard deviation around the average value is fairly high, since each run involves a different dataset and a different random initialisation, hence different value of the likelihood, but the decreasing trend is the same for all of the runs. We also display the average relative square errors on the parameters at the end of the algorithm. They are all small, with reasonably small standard deviation, which indicates that the algorithm consistently recovers correctly the parameters.

To evaluate the impact of the number of Riemann intervals $\varphi(n)$, we run a second experiment where we compare four different profiles over 50 simulations:
\begin{equation*}
 \begin{split}
 \textbf{(low) } \varphi_1(n) &:= n + 1\\
 \textbf{(medium) } \varphi_2(n) &:= n + 100\\
 \textbf{(high) }\varphi_3(n) &:= n + 1000\\
 \textbf{(linear) } \varphi_4(n) &:= 10\times n + 1\, .
 \end{split}
\end{equation*}

\begin{figure}[tbhp]
 \subfloat[]{\includegraphics[width=0.5\linewidth, valign=c]{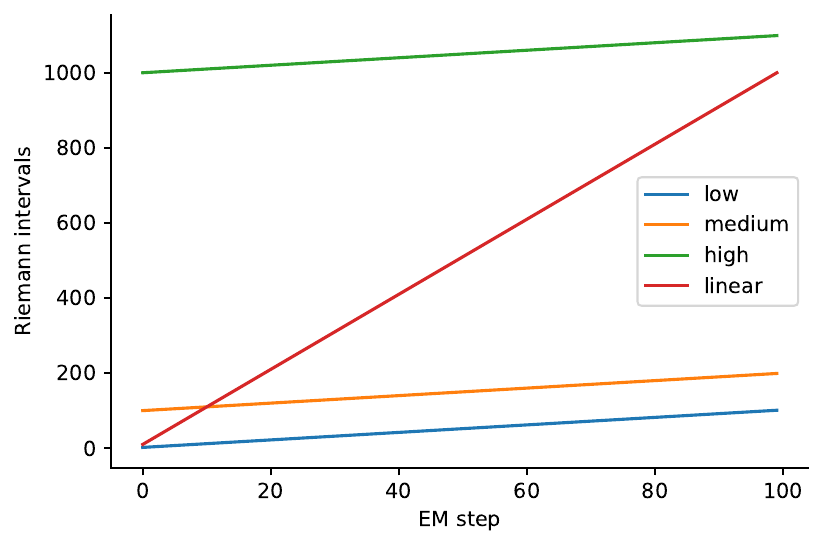}}
 \subfloat[]{\includegraphics[width=0.5\linewidth, valign=c]{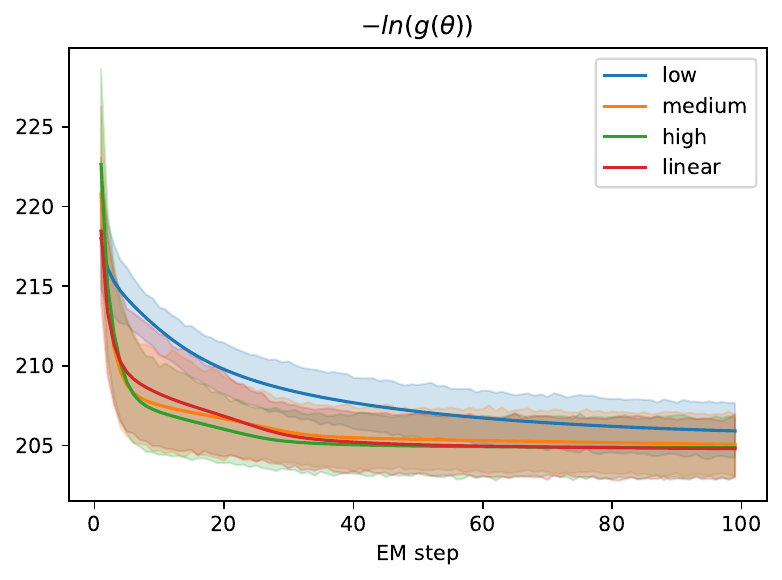}}
\caption{(\textbf{a}). Visual representation of the number of Riemann intervals over the EM steps for each profile $\varphi_i$. The total number of Riemann intervals computed over 100 EM iterations are: 5150 for ``low'', 14,950 for ``medium'', 50,500 for ``linear'' and 104,950 for ``high''. (\textbf{b}). For each profile, average evolution of the negative log-likelihood, with standard deviation, over 50 simulations. The results are fairly similar, in particular between ``medium'', ``high'' and ``linear''.}
\label{fig:riemann_all}
\end{figure}

The results are displayed on \figurename~\ref{fig:riemann_all}. We can see that, despite the very different profiles, the optimisations are very similar. The ``low'' profile performs slightly worst, which indicates that a high number of Riemann intervals is most desirable in practice. As long as this number is high enough, \figurename~\ref{fig:riemann_all} suggests that the performances will not depend too much on the profile.

\vspace{-21pt}

\subsubsection{Application in Two Dimensions}
The difficulty faced by Riemann methods in general is their geometric complexity when the dimension increases. In this section, we propose a similar experiment in two dimensions to show that the method is still functional and practical in that case.

For this 2D-model, we consider two latent independent Beta random variables $z_1 \sim Beta(\alpha_1, 1)$ and $z_2 \sim Beta(\alpha_2, 1)$, and two observed variables defined as: 
\begin{equation*}
\begin{split}
 x_1 &= \lambda_1 z_1 + z_2 + \sigma_1 \epsilon_1 \\
 x_2 &= z_1 + \lambda_2 z_2 + \sigma_2 \epsilon_2 \, ,
\end{split}
\end{equation*}
with $\epsilon_1 \sim \mathcal{N}(0,1)$, $\epsilon_2 \sim \mathcal{N}(0,1)$, and $(z_1, z_2, \epsilon_1, \epsilon_2)$ independent. The 2-dimension version of the Riemann E step with $n$ intervals on each dimension is:
\begin{equation*}
 \tilde{p}_{\theta, n}^{(i)}\parent{\frac{k_1}{n}, \frac{k_2}{n}} = \frac{h^{(i)}(\frac{k_1}{n}, \frac{k_2}{n}; \theta)}{\frac{1}{n^2} \sum_{l_1, l_2=0}^{n-1} h^{(i)}(\frac{l_1}{n}, \frac{l_2}{n}; \theta)} \, .
\end{equation*}

As before, we define an operator $\Psi^{(i)} : \R^{[0, 1]^2} \longrightarrow \R$ such that, for any $f : [0, 1]^2 \longrightarrow \R$:
\begin{equation*}
 \Psi^{(i)} \circ f = \sum_{k_1,k_2=0}^{n-1} \tilde{p}_{\theta, n}^{(i)} \parent{\frac{k_1}{n}, \frac{k_2}{n}} \int_{z_1, z_2=k/n}^{(k+1)/n} f(z_1, z_2) dz \, .
\end{equation*} 

Then, the M step can be expressed as:
\begin{equation*}
 \begin{split}
 \frac{1}{\hat{\alpha}_1} &= -\frac{1}{N} \sum_{i=1}^N \Psi^{(i)} \circ \ln(z_1) \, , \\
 \hat{\lambda}_1 &= \frac{\sum_i \Psi^{(i)} \circ (x_1^{(i)}z_1 - z_2 z_1) }{\sum_i \Psi^{(i)} \circ z_1^2 } \, ,\\
 \hat{\sigma}_1 &=\frac{1}{N}\sum_{i=1}^N \Psi^{(i)} \circ \parent{x_1^{(i)}-\hat{\lambda}_1 z_1 - z_2}^2 \, ,
 \end{split}
\end{equation*}
with symmetric formulas for $\hat{\alpha}_2, \hat{\lambda}_2$ and $\hat{\sigma}_2$. 

For the simulations, we take $(\alpha_1, \alpha_2) = (1, 3)$, $(\lambda_1, \lambda_2) = (10, -10)$ and $(\sigma_1, \sigma_2) = (2, 3)$. From the previous experiment, we keep only the two least costly profiles: ``low'' $\varphi_1(n) := n + 1$ and ``medium'' $\varphi_2(n) := n + 100$. We also consider two new, sub-linear, profiles ``square root'' $\varphi_5(n) := \lfloor \sqrt{n} \rfloor + 1 $ and ``5 square root'' $\varphi_6(n) := 5 \times \lfloor\sqrt{n} \rfloor$. $\varphi_5(n)$ and $\varphi_6(n)$ are designed to have linear complexity even in 2-dimensions.
 
\begin{figure}[tbhp]
 \centering
 \subfloat[]{\includegraphics[width=0.5\linewidth, valign=c]{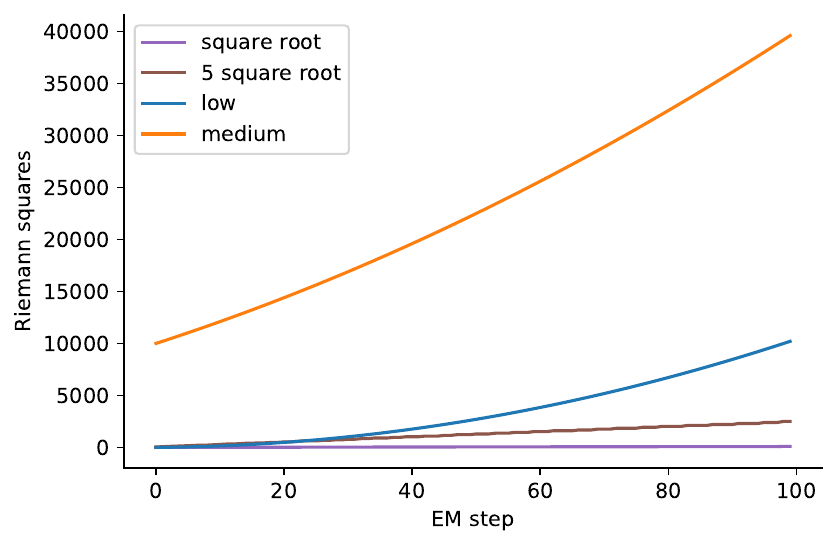}}
 \subfloat[]{\includegraphics[width=0.5\linewidth, valign=c]{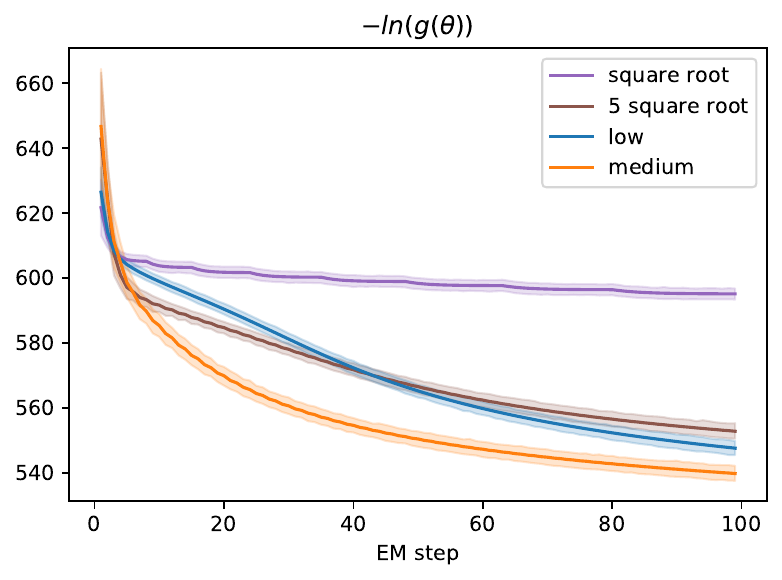}}
\caption{(\textbf{a}). Visual representation of the number of Riemann intervals over the EM steps for each profile $\varphi_i$. In higher dimension, the computational complexity of the profiles are very different. More precisely, the total number of Riemann squares computed over 100 EM iterations are: 4534 for ``square root'', 125,662 for ``5 square root'', 348,550 for ``low'' and 2,318,350 for ``medium''. (\textbf{b}). For each profile, average evolution of the negative log-likelihood, with standard deviation, over 50 simulations. The ``square root'' profile performs poorly. However, the other three are comparable despite their different computational complexities.}
\label{fig:riemann_2D}
\end{figure}

The results of the EM algorithm runs are displayed on \figurename~\ref{fig:riemann_2D}. On the left, we follow the number of Riemann squares mapping the 2D space. The difference in computational complexity between profiles is accentuated by the higher dimension. In particular, ``medium'' performs 6.7 times more computations than ``low'' and 18.4 times more than ``5 square root''. However, on the right of \figurename~\ref{fig:riemann_2D}, we observe that these three profiles perform similar optimisations. This observation justifies cutting computation costs by using lower resolution profiles to compensate the higher dimension. 

\subsection{Tempered EM: Application to Mixtures of Gaussian} \label{sect:exp_tmp_em}

\subsubsection{Context and Experimental Protocol}
In this section, we will assess the capacity of tmp-EM to escape from deceptive local maxima. {We compare the classical EM to tmp-EM with both a monotonous and an oscillating temperature profile} on a very well know toy example: likelihood maximisation within the Gaussian Mixture Model. 
We confront the algorithms to situations where the true classes have increasingly more ambiguous positions, combined with initialisations designed to be hard to escape from. Although the EM is an optimisation procedure, and the log-likelihood reached is a critical metric, in this example, we put more emphasis on the correct positioning of the cluster centroids, that is to say on the recovery of the $\mu_k$. The other usual metrics are also in favour of tmp-EM, and can be found in supplementary materials.

\begin{figure}[tbhp]
 \subfloat[]{\includegraphics[width=0.31\linewidth]{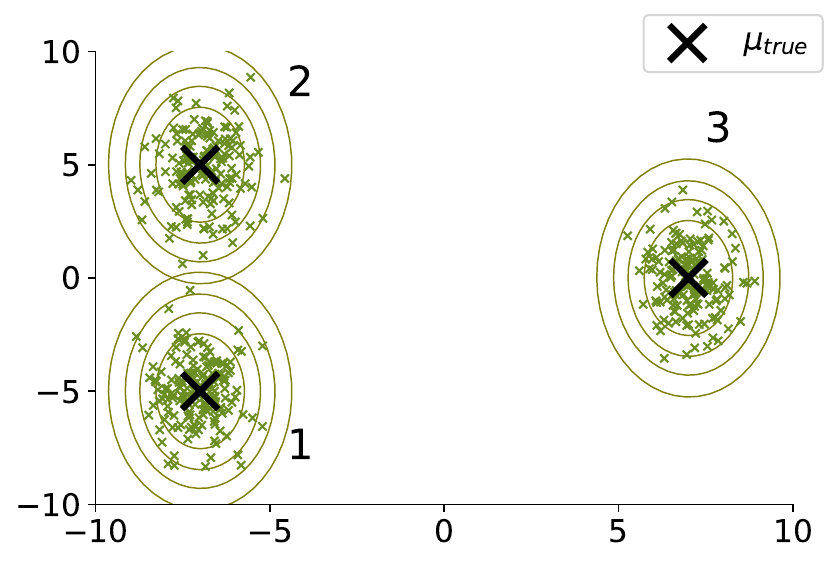}}
 \subfloat[]{\includegraphics[width=0.31\linewidth]{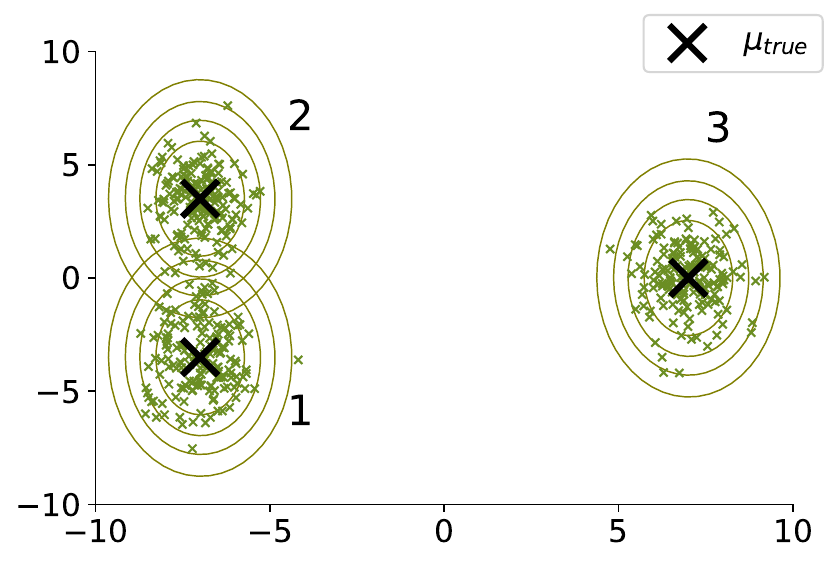}}
 \subfloat[]{\includegraphics[width=0.31\linewidth]{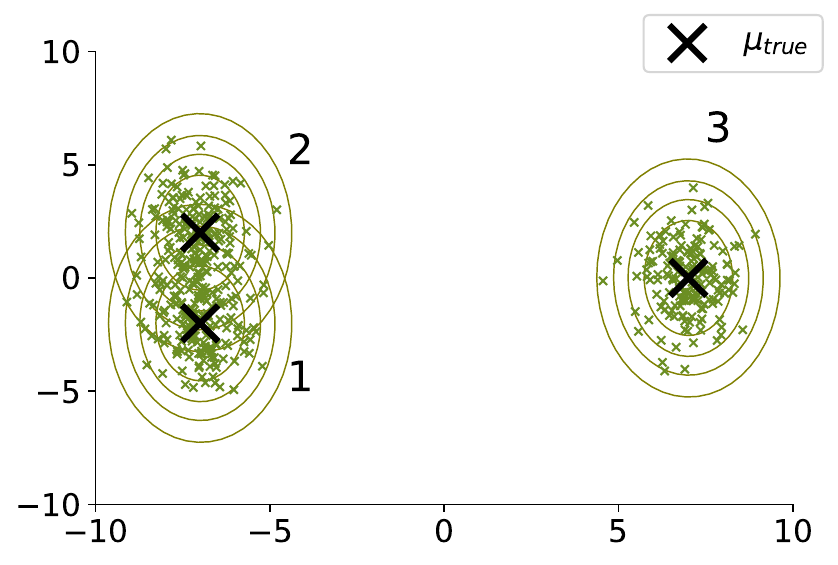}}
\caption{500 sample points from a Mixture of Gaussians with 3 classes. The true centroid of each Gaussian are depicted by black crosses, and their true covariance matrices are represented by the confidence ellipses of level 0.8, 0.99 and 0.999 around the centre. Each sub-figure corresponds to one of the three different versions of the true parameters. From (\textbf{a}) to (\textbf{c}): the true $\mu_k$ of the two left clusters ($\mu_1$ and $\mu_2$) are getting closer while everything else stays identical.}
\label{fig:3_clusters_real_parameters}
\end{figure}

For the sake of comparison, the experimental design is similar to the one in \cite{allassonniere2021new} on the tmp-SAEM. It is as follows: we have three clusters of similar shape and same weight. One is isolated and easily identifiable. The other two are next to one another, in a more ambiguous configuration. \figurename~\ref{fig:3_clusters_real_parameters} represents the three, gradually more ambiguous configurations. Each configuration is called a ``parameter family''.

\begin{figure}[tbhp]
 \subfloat{\includegraphics[width=0.31\linewidth]{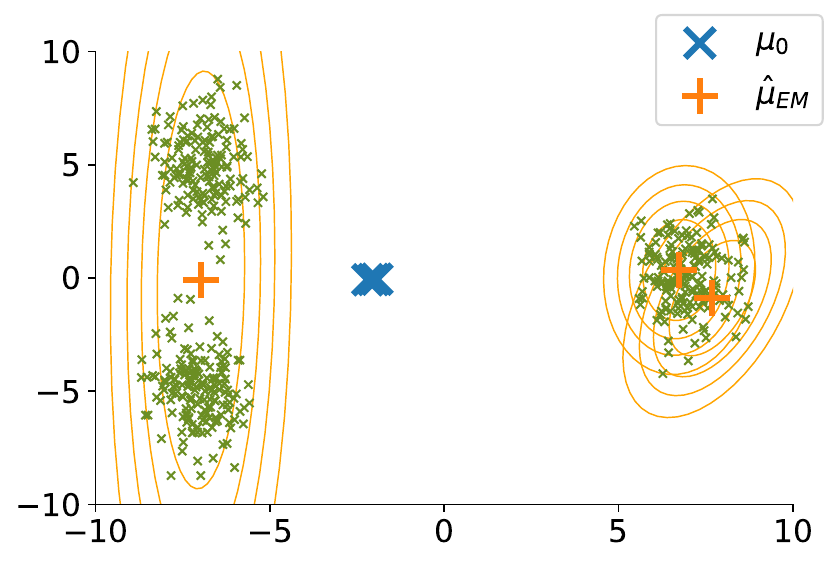}}
 \subfloat{\includegraphics[width=0.31\linewidth]{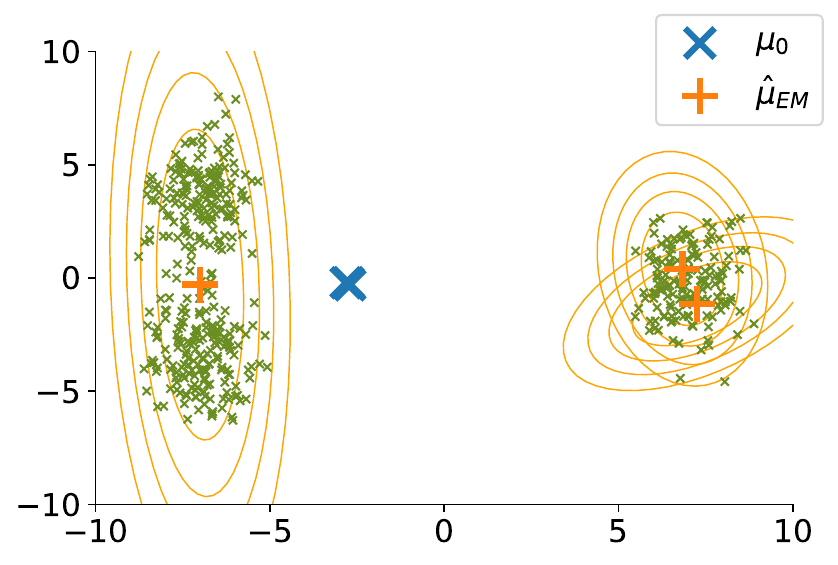}}
 \subfloat{\includegraphics[width=0.31\linewidth]{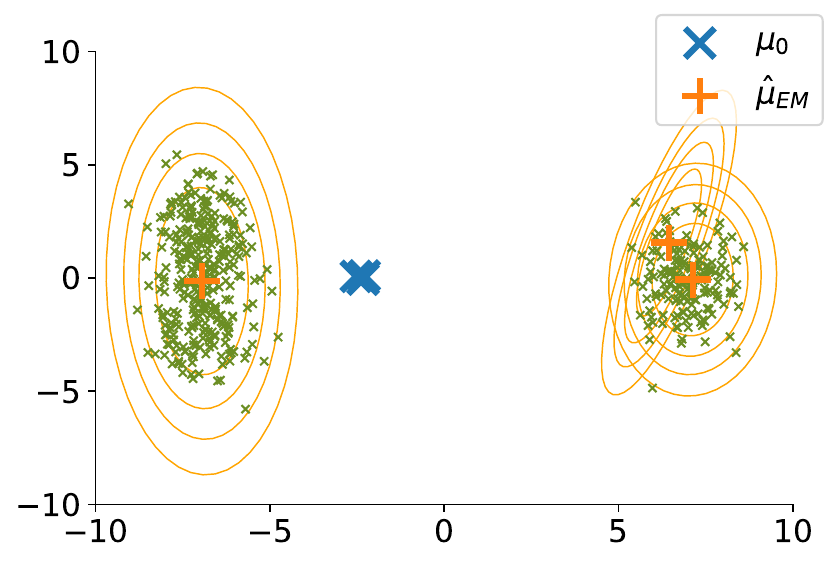}}
 \\
 \subfloat{\includegraphics[width=0.31\linewidth]{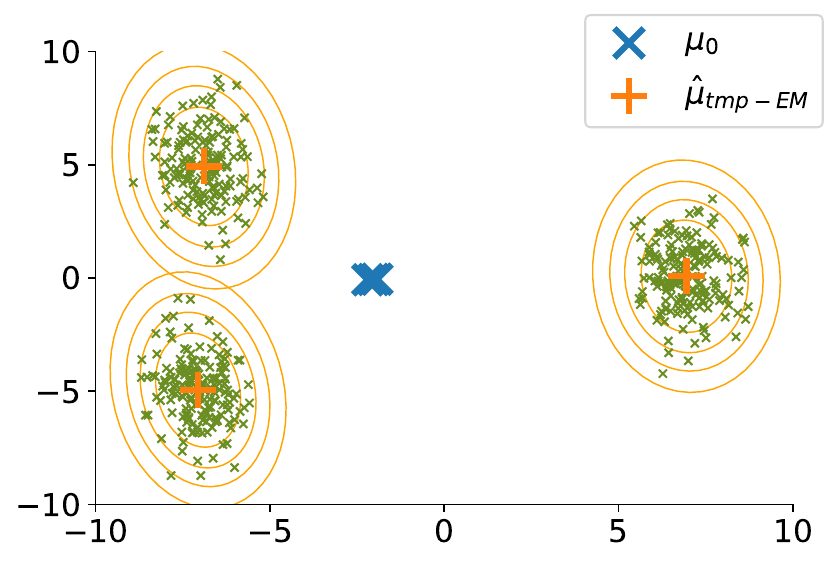}}
 \subfloat{\includegraphics[width=0.31\linewidth]{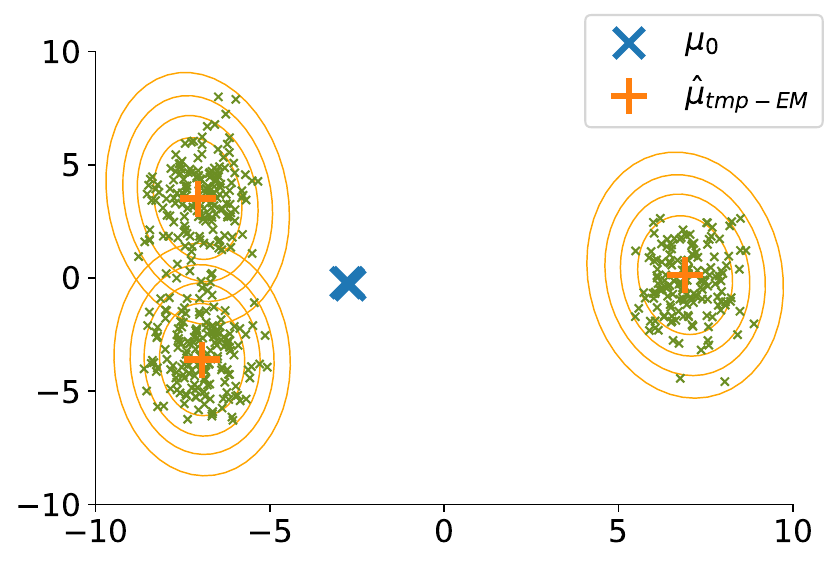}}
 \subfloat{\includegraphics[width=0.31\linewidth]{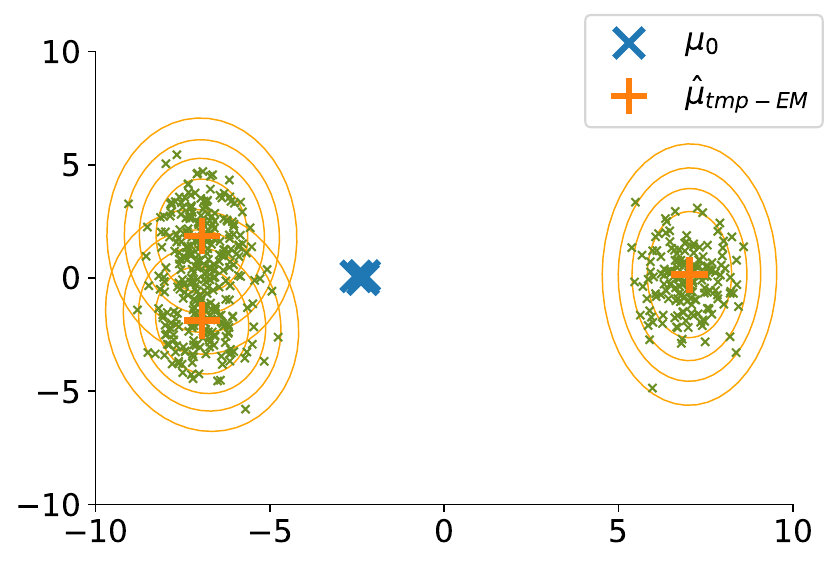}}
\caption{Typical final positioning of the centroids by EM (first row) and tmp-EM {with \textit{oscillating} profile} (second row) \textbf{when the initialisation is made at the barycenter of all data points} (blue crosses). The three columns represent the three gradually more ambiguous parameter sets. Each figure represents the positions of the estimated centroids after convergence of the EM algorithms (orange cross), with their estimated covariance matrices (orange confidence ellipses). In each simulation, 500 sample points were drawn from the real GMM (small green crosses). In those example, tmp-EM managed to correctly identify the position of the three real centroids.}
\label{fig:3_clusters_illustrative_barycenter}
\end{figure}

\begin{figure}[tbhp]
 \subfloat{\includegraphics[width=0.31\linewidth]{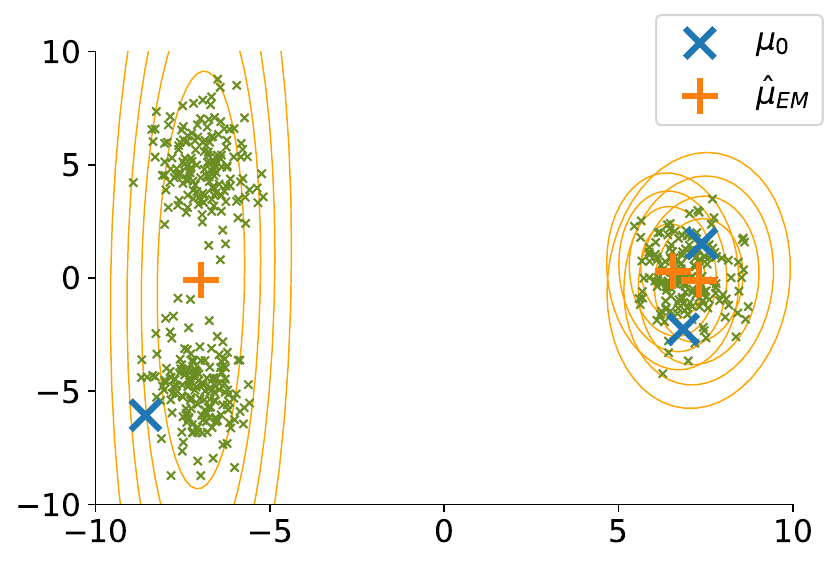}}
 \subfloat{\includegraphics[width=0.31\linewidth]{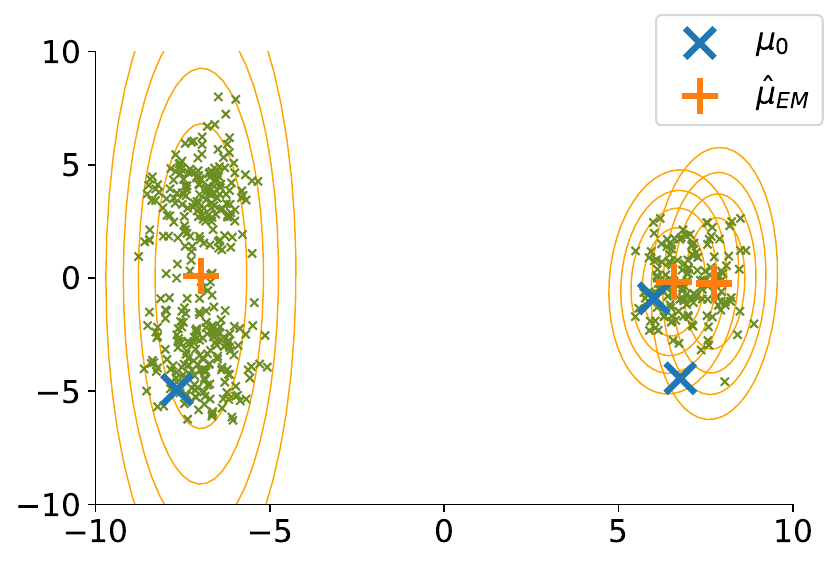}} \subfloat{\includegraphics[width=0.31\linewidth]{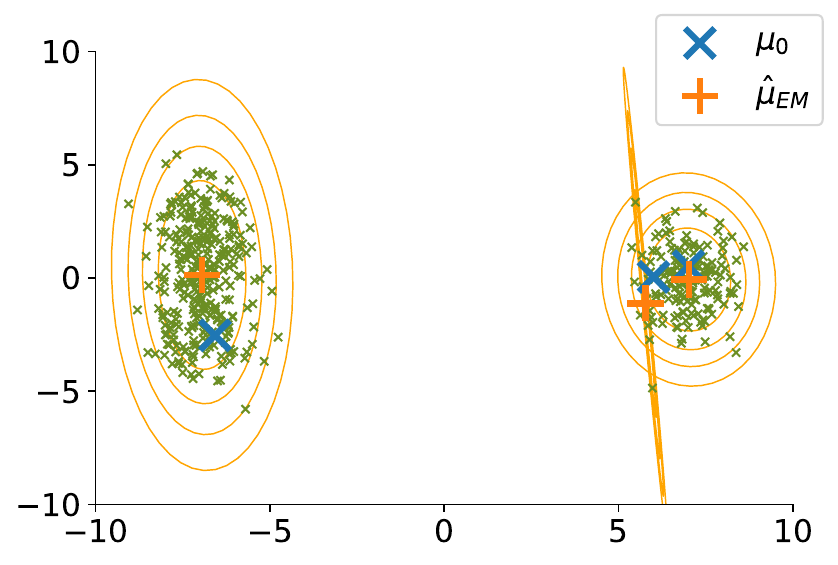}}
 \\
 \subfloat{\includegraphics[width=0.31\linewidth]{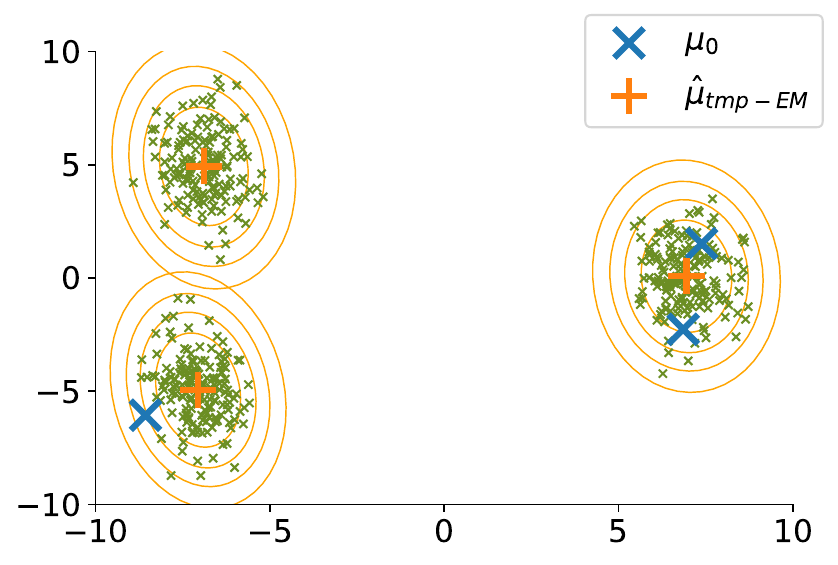}}
 \subfloat{\includegraphics[width=0.31\linewidth]{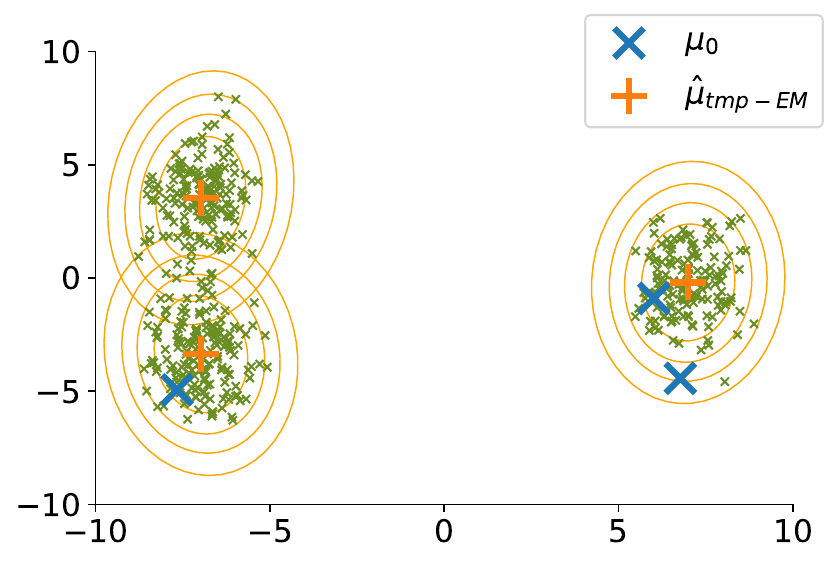}}
 \subfloat{\includegraphics[width=0.31\linewidth]{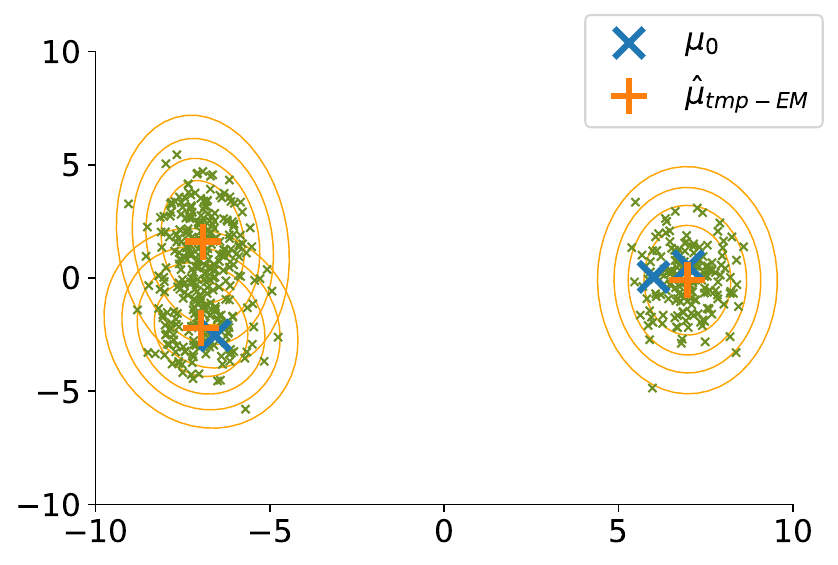}}
\caption{Typical final positioning of the centroids by EM (first row) and tmp-EM {with \textit{oscillating} profile} (second row) \textbf{when the initialisation is made by selecting two points in the isolated cluster and one in the lower ambiguous cluster} (blue crosses). The three columns represent the three gradually more ambiguous parameter sets. Each figure represents the positions of the estimated centroids after convergence of the EM algorithms (orange cross), with their estimated covariance matrices (orange confidence ellipses). In each simulation, 500 sample points were drawn from the real GMM (small green crosses). In those examples, although EM kept two centroids on the isolated cluster, tmp-EM managed to correctly identify the position of the three real centroids.}
\label{fig:3_clusters_illustrative_2v1}
\end{figure}

We use two different initialisation types to reveal the behaviours of the three EMs. The first---which we call ``\textit{barycenter}''---puts all three initial centroids at the centre of mass of all the observed data points. However, none of the EM procedures would move from this initial state if the three GMM centroids were at the exact same position, hence we actually apply a tiny perturbation to make them all slightly distinct. The blue crosses on \figurename~\ref{fig:3_clusters_illustrative_barycenter} represent a typical \textit{barycenter} initialisation. With this initialisation method, we assess whether the EM procedures are able to correctly estimate the positions of the three clusters, despite the ambiguity, when starting from a fairly neutral position, providing neither direction nor misdirection. On the other hand, the second initialisation type - which we call ``\textit{2v1}'' - is voluntarily misguiding the algorithm by positioning two centroids on the isolated right cluster and only one centroid on the side of the two ambiguous left clusters. The blue crosses on \figurename~\ref{fig:3_clusters_illustrative_2v1} represent a typical \textit{2v1} initialisation. This initialisation is intended to assess whether the methods are able to escape the potential well in which they start and make theirs centroids traverse the empty space between the left and right clusters to reach their rightful position. For each of the three parameter families represented on \figurename~\ref{fig:3_clusters_real_parameters}, 1000 datasets with 500 observations each are simulated, and the three EMs are ran with both the \textit{barycenter} and the \textit{2v1} initialisation.

{For tmp-EM, we try two profiles. First, a simple \textit{decreasing} exponential profile as seen in \cite{ueda1998deterministic}: $T_n = 1 + (T_0-1) \, \exp(-r.n)$. Through a grid search, the values \mbox{$T_0=5$}, $ r=2$ for the \textit{barycenter} initialisation and $T_0=100\, r=1.5$ for the \textit{2v1} initialisation are picked for this profile. Since Theorem \ref{thm:tempered_EM} only requires $T_n \longrightarrow 1$, we also propose an \textit{oscillating} profile inspired from \cite{allassonniere2021new}. The exact formula of these oscillations is:} $T_n = th(\frac{n}{2r}) + (T_0 - b \frac{2 \sqrt{2}}{3 \pi})\, a^{n/r} + b \, sinc(\frac{3 \pi}{4}+ \frac{n}{r})$. Where $0< T_0$, $0 < r$, $0 < b$ and $ 0<a<1$. The oscillations in this profile are meant to achieve a two-regimes behaviour. When the temperature reaches low values, the convergence speed is momentarily increased which has the effect of ``locking-in'' some of the most obviously good decisions of the algorithm. Then, the temperature is re-increased to continue the exploration on the other, more ambiguous, parameters. Those two regimes are alternated in succession with gradually smaller oscillations, resulting in a multi-scale procedure that ``locks-in'' gradually harder decisions. For some hyper-parameter combinations, the sequence $T_n$ can have a (usually small) finite number of negative values. Since only the asymptotic behaviour of $T_n$ is the step $n$ matters for convergence, then the theory allows a finite number of negative values. However, in practice, at least for the sake of interpretation, one may prefer to use only positive values for $T_n$. In which case, one can either restrain themselves to parameter combinations that result in no negatives values for $T_n$, or enforce positivity by taking $T_n \leftarrow max(T_n, \epsilon)$ with a certain $\epsilon>0$. 

For our experiments, we select the hyper-parameter values with a grid-search. The ``normalised'' $sinc$ function is used $sinc(x)=sin(\pi x)/(\pi x)$ and the chosen tempering parameters are $T_0=5,\, r=2,\, a=0.6,\, b=20$ for the experiments with the \textit{barycenter} initialisation, and $T_0=100,\, r=1.5,\, a=0.02,\, b=20$ for the \textit{2v1} initialisation. We have two different sets of tempering hyper-parameters values, one for each of the two very different initialisation types. However, these values then remain the same for the three different parameter families and for every data generation within them. This underlines that the method is not excessively sensitive to the tempering parameters, and that the prior search for good hyper-parameter values is a worthwhile time investment. Likewise, a simple experiment with 6 clusters, in supplementary materials, demonstrates that the same hyper-parameters can be kept over different initialisation (and different data generations as well) when they were made in a non-adversarial way, by drawing random initial centroids uniformly among the data points.

\vspace{-15pt}
 
\subsubsection{Experimental Results Analysis}
{In this section, we analyse the results of EM, \textit{decreasing} tmp-EM and \textit{oscillating} tmp-EM over all the simulations.}

{First, an illustration: Figures~\ref{fig:3_clusters_illustrative_barycenter} and \ref{fig:3_clusters_illustrative_2v1} depict the final states of EM and \textit{oscillating} tmp-EM} on one typical simulation for each of the three ambiguity level (the three parameter families) starting from the \textit{barycenter} and \textit{2v1} initialisation respectively. The simulated data are represented by the green crosses. The initial centroids are in blue. The orange cross represents the estimated centroids positions $\hat{\mu}_k$, and the orange confidence ellipses are visual representations of the estimated covariance matrices $\hat{\Sigma}_k$. In supplementary materials, we show step by step the path taken by the estimated parameters of tmp-EM before convergence, providing much more detail on the method's behaviours. {These illustrative examples show \textit{oscillating} tmp-EM better succeeding at clustering recovery than the classical EM. The results over all simulations are aggregated in Table \ref{tab:3_clusters_quantitative}, and confirm this observation.}

Table \ref{tab:3_clusters_quantitative} presents the average and the standard deviation of the relative $l_2$ error on $\mu_k$ of the EMs. For each category, the better result over {the three EM }is highlighted in bold. The recovery of the true class averages $\mu_k$ is spotlighted as it is the essential success metric for this experiment. {More specifically, class 1 and 2, the two leftmost classes, are the hardest to correctly recover and the ones whose estimation is the differentiating criterion between the algorithms. Indeed, as can be seen in Table \ref{tab:3_clusters_quantitative}, $\mu_3$ is always well estimated by all methods. Hence, in the following, we discuss the error on $\mu_1$ and $\mu_2$.} 

{With the \textit{barycenter} initialisation, classical EM and \textit{decreasing} tmp-EM have similar average relative error levels. With classical EM actually being slightly better. However, \textit{oscillating} tmp-EM is much better than both of them, with error levels smaller by a factor of 10 on parameter families 1 and 2, and by a factor of 5 on parameter family 3. The standard deviation of \textit{oscillating} tmp-EM is also lower, by a factor of roughly 3 on parameter families 1 and 2, and by a factor of 2 on parameter family 3. With the \textit{2v1} initialisation, all error levels are higher. This time, \textit{decreasing} tmp-EM is better in average than classical EM by a factor of 1.7 to 1.4, depending on the parameter family. In turn, \textit{oscillating} tmp-EM is better than \textit{decreasing} tmp-EM by a factor of 3.1 to 3.4 depending on the parameter family. Its standard deviation is also lower by a factor of about 2.

Overall, \textit{oscillating} tmp-EM dominates the simulation. Its error rates on the recovery of $\mu_1$ and $\mu_2$ are always the best, and they remain at low levels even with the most adversarial initialisations. To bolster this last point, we underline that the highest relative error reached by \textit{oscillating} tmp-EM over all the various scenarios (0.39 on parameter family 3 with \textit{2v1} initialisation) is still lower than the lowest relative error of both classical EM (0.52 on parameter family 1 with \textit{barycenter} initialisation) and \textit{decreasing} tmp-EM (0.60 on parameter family 1 with \textit{barycenter} initialisation). 
}

\begin{table*}[tbhp]
\begin{center}
{\footnotesize
\caption{Average and standard deviation of the relative error on $\mu_k$, $\frac{\norm{\hat{\mu}_k - \mu_k}^2}{\norm{\mu_k}^2}$, made by EM, {tmp-EM with \textit{decreasing} temperature and tmp-EM with \textit{oscillating} temperature} over 1000 simulated dataset with two different initialisations. The three different parameter families, described in \figurename~\ref{fig:3_clusters_real_parameters}, correspond to increasingly ambiguous positions of classes 1 and 2. For both initialisations type, the identification of these two clusters is drastically improved by the tempering. {Best results highlighted in \textbf{bold}.}} \label{tab:3_clusters_quantitative}
\begin{tabular}{lccccccc}
\toprule
 & & \multicolumn{2}{c}{\textbf{EM}}& \multicolumn{2}{c}{\textbf{tmp-EM (\textit{Decreasing}} \boldmath{$T$}\textbf{)}} & \multicolumn{2}{c}{\textbf{tmp-EM (\textit{Decreasing
oscillating}} \boldmath{$T$}\textbf{)}}\\
\midrule
Parameter\\ Family & cl. & barycenter & 2v1 & barycenter& 2v1 & barycenter& 2v1 \\
\midrule
\multirow{3}{*}{1} 
& 1 & 0.52 (1.01) & 1.52 (1.24) & 0.60	(1.08) & 0.87 (1.20) &\textbf{0.04 (0.26)} & \textbf{0.29 (0.64)}\\
& 2 & 0.55 (1.05) & 1.53 (1.25) & 0.64	(1.10) & 0.96 (1.25) & \textbf{0.05 (0.31)} & \textbf{0.30 (0.64) } \\
& 3 & \textbf{0.01 (0.06)} &0.01 (0.03) & 0.01	(0.10) & \textbf{0.01 (0.02)} & 0.03 (0.17)& 0.03 (0.19) \\
\midrule
\multirow{3}{*}{2}
& 1 & 1.00 (1.42)& 1.69 (1.51) & 0.96 (1.41) & 1.10 (1.46) &\textbf{0.09 (0.47)}& \textbf{0.37 (0.86)} \\
& 2 & 1.03 (1.44)& 1.71 (1.52) & 1.08	(1.46) & 1.11 (1.46) &\textbf{0.12 (0.57)}& \textbf{0.32 (0.79)} \\
& 3 & 0.01 (0.05)& \textbf{0.02 (0.03)} & 0.01 \textbf{(0.03)} & 0.01 (0.04) & \textbf{0.01} (0.05) & 0.04 (0.22) \\
\midrule
\multirow{3}{*}{3}
& 1 & 1.56 (1.75) & 1.79 (1.77) & 1.63	(1.76) & 1.38 (1.71) &\textbf{0.31 (0.97)}& \textbf{0.39 (0.98)}\\
& 2 & 1.51 (1.74) & 1.88 (1.76) & 1.52	(1.74) & 1.30 (1.68) &\textbf{0.30 (0.93)}& \textbf{0.39 (0.97)}\\
& 3 & 0.02 (0.04) & \textbf{0.02 (0.04)} & 0.01 \textbf{(0.03)}	& 0.02 (0.06) & \textbf{0.01} (0.04)& 0.07 (0.30)\\
\bottomrule
\end{tabular}
}
\end{center}
\end{table*}

\subsection{Tempered Riemann Approximation EM: Application to a Gaussian Model with Beta Prior} \label{sect:exp_tmp_riemann_em}
We illustrate the method with the model of Section \ref{sect:riemann_em_example}:
\begin{equation*}
 h(z; \theta) =\frac{\alpha z^{\alpha-1}}{\sqrt{2 \pi \sigma^2}} exp\parent{-\frac{(y-\lambda z)^2}{2 \sigma^2}}\, .
\end{equation*}

We apply the tempered Riemann approximation. As in Section \ref{sect:riemann_em_example}, the resulting conditional probability density is a step function defined by the $n$ different values it takes on $[0,1]$. For the observation $x^{(i)}$, $\forall k \in \llbracket{0, n-1} \rrbracket$:
\begin{equation*}
 \tilde{p}_{\theta, n}^{(i)}\parent{\frac{k}{n}} = \frac{h^{(i)}\parent{\frac{k}{n}; \theta}^{\frac{1}{T_n}}}{\frac{1}{n} \sum_{l=0}^{n-1} h^{(i)}\parent{\frac{l}{n}; \theta}^{\frac{1}{T_n}}} \, .
\end{equation*}

\begin{figure}[tbhp]
\centering 
 \includegraphics[width=\linewidth]{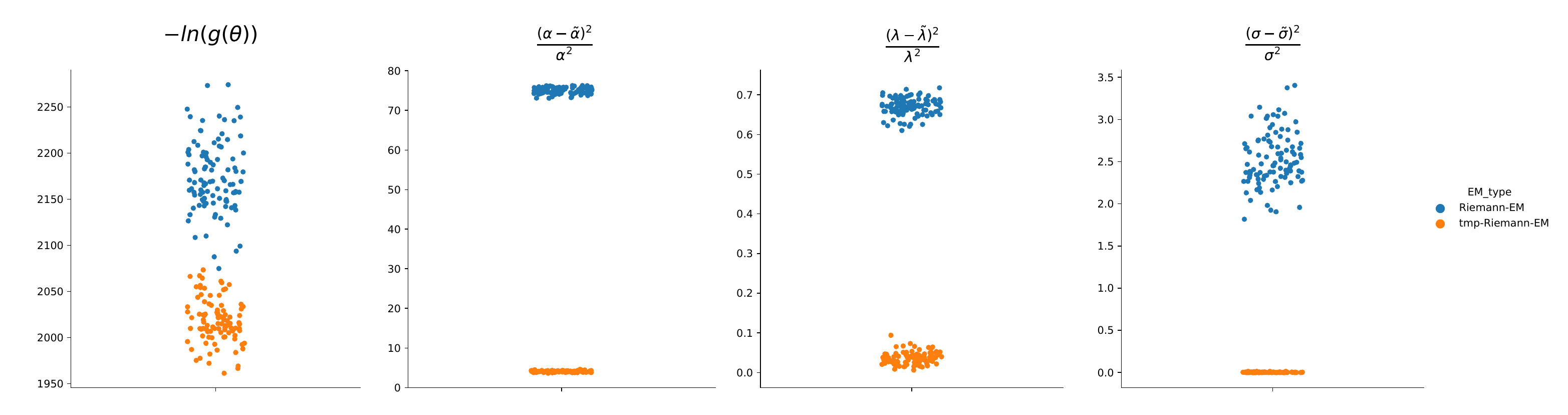}
\caption{Results over many simulations of the Riemann EM and tmp-Riemann EM on the Beta-Gaussian model. The tempered Riemann EM reaches relative errors on the real parameters that are orders of magnitude below the Riemann EM with no temperature. The likelihood reached is also lower with the tempering.}
\label{fig:tmp_riemann}
\end{figure}

The M step, seen in Equation~\eqref{eq:gaussian_beta_M_step}, is unchanged. We compare the tempered Riemann EM to the simple Riemann EM on a case where the parameters are ambiguous. With real parameters $\alpha = 0.1, \lambda =10, \sigma = 0.8$, for each of the 100 simulations, the algorithms are initialised at $\alpha_0 = 10, \lambda_0 =1, \sigma_0 = 7$. The initialisation is somewhat adversarial, since the mean and variance of the marginal distribution of $y$ are approximately the same with the real of the initialisation parameter, even though the distribution is different. \figurename~\ref{fig:tmp_riemann} shows that the tempered Riemann EM better escapes the initialisation than the regular Riemann EM, and reaches errors on the parameters orders of magnitude below. The tempering parameters are here $T_0 =150, r =3, a=0.02, b=40$. 
\vspace{-9pt}

\section{Discussion and Conclusions}\label{sect:conclusions}
We proposed the Deterministic Approximate EM class to bring together the many possible deterministic approximations of the E step. We proved a unified Theorem, with mild conditions on the approximation, which ensures the convergence of the algorithms in this class. Then, we showcased members of this class that solve the usual practical issues of the EM algorithm. For intractable E steps {in low dimension}, we introduced the Riemann approximation EM, a less parametric and deterministic alternative to the extensive family of MC-EM. We showed on an empirical intractable example how the Riemann approximation EM was able to increase the likelihood and recover every parameter in a satisfactory manner with its simplest design, and no hyper parameter optimisation.

{Second, we studied the tempered EM, introduced by \citet{ueda1998deterministic} to escape the attractive sub-optimal local extrema of non-convex objectives. We proved that tmp-EM is a specific case of the Deterministic Approximate EM, benefiting from the convergence property as long as the temperature profile converges towards 1. This mild condition justifies the use of many more temperature profiles than the ones tried in \cite{ueda1998deterministic,naim2012convergence}. To illustrate the interest of complex, non-monotonous, temperature profiles, we demonstrated on experiments with adversarial initial positions the superiority of an \textit{oscillating} profile over a simple \textit{decreasing}~one.}

Finally, we added the Riemann approximation in order to apply the tempering in intractable cases. We were then able to show that the tmp-Riemann approximation massively improved the performances of the Riemann approximation, when the initialisation is ambiguous.

Future works will improve both methods. The Riemann approximation will be generalised to be applicable even when the latent variable is not bounded, and an intelligent slicing of the integration space will improve the computational performances in high dimension. Regarding the tempered EM, since the theory allows the usage of any temperature profile, the natural next step is to look for efficient profiles with few hyper-parameters for fast tuning. Afterwards, implementing an adaptive tuning of the temperature parameters during the procedure will remove the necessity for preliminary grid search altogether.

\vspace{6pt} 




\paragraph{Acknowledgements}{The research leading to these results has received funding from the European Research Council (ERC) under grant agreement No 678304, European Union’s Horizon 2020 research and innovation program under grant agreement No 666992 (EuroPOND) and No 826421 (TVB-Cloud), and the French government under management of Agence Nationale de la Recherche as part of the ``Investissements d'avenir'' program, reference ANR-19-P3IA-0001 (PRAIRIE 3IA Institute) and reference ANR-10-IAIHU-06 (IHU-A-ICM).}

\appendix

\section{Proofs of the Two Main Theorems} \label{sect:proofs}
In this Appendix, we provide in full details the proofs that were sketched in the main body of the paper. 
Section \ref{sect:proofs_main} details the proof of convergence for the Deterministic Approximate EM algorithm, our central result, Theorem \ref{thm:main} of the paper. Section \ref{sect:proofs_tempered} details the proof of convergence for tmp-EM, Theorem \ref{thm:tempered_EM} of the paper. 

\subsection{Proof of the General Theorem} \label{sect:proofs_main} 
In this Section, we prove Theorem \ref{thm:main}, which guarantees the convergence of the Deterministic Approximate EM algorithm.

The proof is based on the application of Propositions \ref{prop:9} and \ref{prop:11}, taken from \citet{fort2003convergence}.

We need to prove that, under the conditions of Theorem \ref{thm:main}, we verify the conditions of Proposition Proposition \ref{prop:9} and Proposition \ref{prop:11}. Then we will have the results announced in Theorem \ref{thm:main}.

\subsubsection{Verifying the Conditions of Proposition \ref{prop:11}}
$g$ is the likelihood function of a model of the curved exponential family. Let $T$ be the point to point map describing the transition between $\theta_n$ and $\theta_{n+1}$ in the exact EM algorithm. $\L$ the set of stationary points by $T$: $\mathcal{L} := \brace{\theta \in \Theta | T(\theta) = \theta }$. (Note that if $g$ is differentiable, the general properties of the EM tell us that its critical points of $g$ are the stationary points: $\mathcal{L} = \brace{\theta \in \Theta | \nabla g (\theta) = 0 }$). Additionally, $g$ is a $C^0$ Lyapunov function associated to $(T, \L)$. Let $\brace{\theta_n}_n$ be the sequence defined by the stable approximate EM with $\{F_n\}_{n \in \N}$ our sequence of point to point maps.

We verify that under this framework---and with the assumptions of Theorem \ref{thm:main}---we check the conditions of Proposition \ref{prop:11}.

As in \cite{fort2003convergence}, \emph{M1--3} directly implies \emph{C1--2}. 

Let us show that we have the last two conditions for Proposition \ref{prop:11}: 
\begin{equation} \label{eq:prop_11_condition_1}
 {\forall \theta \in K_0, \quad \underset{n \to \infty}{\lim}\,|g \circ F_n - g \circ T|(\theta) = 0} \, ,
\end{equation}
and
\begin{equation} \label{eq:prop_11_condition_2}
 {\forall \text{ compact } K \subseteq \Theta, \quad \underset{n \to \infty}{\lim} |g \circ F_n(\theta_n) - g \circ T (\theta_n)| \mathds{1}_{\theta_n \in K} = 0\, .}
\end{equation}

We focus on \eqref{eq:prop_11_condition_2}, since \eqref{eq:prop_11_condition_1} is easier to verify and will come from the same reasoning. The first steps are similar to \cite{fort2003convergence}. We underline the most consequent deviations from the proof of \cite{fort2003convergence} when they occur.

\paragraph{Equivalent Formulation of the Convergence} 
We write Equation~\eqref{eq:prop_11_condition_2} under an equivalent form. Let $\Tilde{S}_n(\theta) := \brace{\int_z S_u(z) \Tilde{p}_{\theta, n}(z) dz}_{u=1}^q$ and $\bar{S}(\theta) := \brace{\int_z S_u(z) p_{\theta}(z) dz}_{u=1}^q$. Then $F_n(\theta_n) = \Hat{\theta}(\Tilde{S}_n(\theta_n))$ and $T(\theta_n) = \Hat{\theta}(\bar{S}(\theta_n))$. Hence $|g \circ F_n(u_n) - g \circ T (u_n)| = |g \circ \Hat{\theta}(\Tilde{S}_n(\theta_n)) - g \circ \Hat{\theta}(\bar{S}(\theta_n))|$. To show Equation~\eqref{eq:prop_11_condition_2}:
\begin{equation*}
 {|g \circ \Hat{\theta}(\Tilde{S}_n(\theta_n)) - g \circ \Hat{\theta}(\bar{S}(\theta_n)) | \mathds{1}_{\theta_n \in K} \underset{n \to \infty}{\to} 0\, ,}
\end{equation*}
it is sufficient and necessary to have:
\begin{equation*}
 \forall \epsilon > 0, \, \exists N \in \N, \, \forall n \geq N, \, |g \circ \Hat{\theta}(\Tilde{S}_n(\theta_n)) - g \circ \Hat{\theta}(\bar{S}(\theta_n)) | \mathds{1}_{\theta_n \in K} \leq \epsilon \, .
\end{equation*}

An other equivalent formulation is that there are a finite number of integers $n$ such that $|g \circ \Hat{\theta}(\Tilde{S}_n(\theta_n)) - g \circ \Hat{\theta}(\bar{S}(\theta_n)) | \mathds{1}_{\theta_n \in K} > \epsilon$, in other words:
\begin{equation*}
 \forall \epsilon > 0, \, \sum_{n=1}^{\infty} \mathds{1}_{|g \circ \Hat{\theta}(\Tilde{S}_n(\theta_n)) - g \circ \Hat{\theta}(\bar{S}(\theta_n)) | \mathds{1}_{\theta_n \in K} > \epsilon} < \infty \, .
\end{equation*}
\paragraph{Use the Uniform Continuity} We aim to relate the proximity between the images $g \circ \Hat{\theta}$ of to the proximity between the antecedents of $g \circ \Hat{\theta}$. The function $g \circ \Hat{\theta} : \R^q \to \R$ is continuous, but not necessarily uniformly continuous on $\R^q$. As a consequence, we will need to restrict ourselves to a compact to get uniform continuity properties. We already have a providen compact $K$. $\Tilde{S} : \Theta \to \R^l$ is continuous, hence $S(K)$ is a compact as well. Let $\delta$ be a strictly positive real number. Let $\bar{S}(K, \delta) := \brace{s \in \R^q \Bigg| \underset{t \in K}{inf} \, ||\bar{S}(t)-s|| \leq \delta}$. Where we use any norm $||.||$ on $\R^q$ since they are all equivalent. $\bar{S}(K, \delta)$ is a compact set as well. As a consequence $g \circ \theta$ is uniformly continuous on $\bar{S}(K, \delta)$, which means that:

\begin{equation} \label{eq:uniform_continuity}
 \forall \epsilon > 0, \, \exists \eta(\epsilon, \delta) > 0, \, \forall x, y \in \bar{S}(K, \delta), \, \norm{x-y} \leq \eta(\epsilon, \delta) \implies | g \circ \hat{\theta}(x) - g \circ \hat{\theta}(y) | \leq \epsilon \, .
\end{equation}

Let us show that, with $\alpha := min(\delta, \eta(\epsilon, \delta))$, $\forall n,$
\begin{equation} \label{eq:new_implication}
 |g \circ \Hat{\theta}(\Tilde{S}_n(\theta_n)) - g \circ \Hat{\theta}(\bar{S}(\theta_n)) | \mathds{1}_{\theta_n \in K} > \epsilon \implies \norm{\Tilde{S}_n(\theta_n) - \bar{S}(\theta_n) } \mathds{1}_{\theta_n \in K} > \alpha \, .
\end{equation}

To that end, we show that:
\begin{equation*}
 \norm{\Tilde{S}_n(\theta_n) - \bar{S}(\theta_n) } \mathds{1}_{\theta_n \in K} \leq \alpha \implies |g \circ \Hat{\theta}(\Tilde{S}_n(\theta_n)) - g \circ \Hat{\theta}(\bar{S}(\theta_n)) | \mathds{1}_{\theta_n \in K} \leq \epsilon \, .
\end{equation*}

Let us assume that $\norm{\Tilde{S}_n(\theta_n) - \bar{S}(\theta_n) } \mathds{1}_{\theta_n \in K} \leq \alpha$.

If $\theta_n \notin K$, then $|g \circ \Hat{\theta}(\Tilde{S}_n(\theta_n)) - g \circ \Hat{\theta}(\bar{S}(\theta_n)) | \mathds{1}_{\theta_n \in K} = 0 \leq \epsilon$.

If, in contrary, $\theta_n \in K$, then $\bar{S}(\theta_n) \in \bar{S}(K) \subset \bar{S}(K, \delta)$.

Since $\norm{\Tilde{S}_n(\theta_n) - \bar{S}(\theta_n) } = \norm{\Tilde{S}_n(\theta_n) - \bar{S}(\theta_n) } \mathds{1}_{\theta_n \in K} \leq \alpha \leq \delta$, then $\Tilde{S}_n(\theta_n) \in \bar{S}(K, \delta)$.\\ Since $(\bar{S}(\theta_n), \Tilde{S}_n(\theta_n)) \in \bar{S}(K, \delta)^2$ and $\norm{\Tilde{S}_n(\theta_n) - \bar{S}(\theta_n) } \leq \alpha \leq \eta(\epsilon, \delta)$, then we get from Equation~\eqref{eq:uniform_continuity}
\begin{equation*}
 |g \circ \Hat{\theta}(\Tilde{S}_n(\theta_n)) - g \circ \Hat{\theta}(\bar{S}(\theta_n)) | \mathds{1}_{\theta_n \in K} \leq \epsilon \, .
\end{equation*}

In both cases, we get that:
\begin{equation*}
 \norm{\Tilde{S}_n(\theta_n) - \bar{S}(\theta_n) } \mathds{1}_{\theta_n \in K} \leq \alpha \implies |g \circ \Hat{\theta}(\Tilde{S}_n(\theta_n)) - g \circ \Hat{\theta}(\bar{S}(\theta_n)) | \mathds{1}_{\theta_n \in K} \leq \epsilon \, ,
\end{equation*}
which proves Equation~\eqref{eq:new_implication}.

\paragraph{Sufficient Condition for Convergence} 
We use Equation~\eqref{eq:new_implication} to find a sufficient condition for \eqref{eq:prop_11_condition_2}. This part differs from \cite{fort2003convergence} as our approximation is not defined as a random sum. Equation~\eqref{eq:new_implication} is equivalent to
\begin{equation*}
 \mathds{1}_{|g \circ \Hat{\theta}(\Tilde{S}_n(\theta_n)) - g \circ \Hat{\theta}(\bar{S}(\theta_n)) | \mathds{1}_{\theta_n \in K} > \epsilon} \leq \mathds{1}_{\norm{\Tilde{S}_n(\theta_n) - \bar{S}(\theta_n) } \mathds{1}_{\theta_n \in K} > \alpha} \, .
\end{equation*}

From that, we get 
\begin{equation*}
 \forall \epsilon > 0,\, \exists \alpha > 0 \, \sum_{n=1}^{\infty} \mathds{1}_{|g \circ \Hat{\theta}(\Tilde{S}_n(\theta_n)) - g \circ \Hat{\theta}(\bar{S}(\theta_n)) | \mathds{1}_{\theta_n \in K} > \epsilon} \leq \sum_{n=1}^{\infty} \mathds{1}_{\norm{\Tilde{S}_n(\theta_n) - \bar{S}(\theta_n) } \mathds{1}_{\theta_n \in K} > \alpha} \, .
\end{equation*}

As a consequence, if 
\begin{equation*}
 \forall \alpha > 0, \, \sum_{n=1}^{\infty} \mathds{1}_{\norm{\Tilde{S}_n(\theta_n) - \bar{S}(\theta_n) } \mathds{1}_{\theta_n \in K} > \alpha} < \infty
\end{equation*}

Then
\begin{equation*}
 \forall \epsilon > 0,\, \sum_{n=1}^{\infty} \mathds{1}_{|g \circ \Hat{\theta}(\Tilde{S}_n(\theta_n)) - g \circ \Hat{\theta}(\bar{S}(\theta_n)) | \mathds{1}_{\theta_n \in K} > \epsilon} < \infty
\end{equation*}

In other, equivalent, words:
\begin{equation} \label{eq:proof_sufficient_condition}
\begin{split}
 \text{If}& \quad \norm{\Tilde{S}_n(\theta_n) - \bar{S}(\theta_n) } \mathds{1}_{\theta_n \in K} \underset{n \to \infty}{\longrightarrow} 0 \\
 \text{Then}& \quad {|g \circ \Hat{\theta}(\Tilde{S}_n(\theta_n)) - g \circ \Hat{\theta}(\bar{S}(\theta_n)) | \mathds{1}_{\theta_n \in K} \underset{n \to \infty}{\longrightarrow} 0} \, .
\end{split}
\end{equation}

Hence, having for all compact sets $K \subset \Theta, \; \norm{\Tilde{S}_n(\theta_n) - \bar{S}(\theta_n) } \mathds{1}_{\theta_n \in K} \underset{n \to \infty}{\longrightarrow} 0$ is sufficient to have the desired condition \eqref{eq:prop_11_condition_2}. Similarly, we find that $\forall \theta \in K_0$:
\begin{equation} \label{eq:proof_sufficient_condition_2}
 \begin{split}
 &\norm{\Tilde{S}_n(\theta) - \bar{S}(\theta) } \underset{n \to \infty}{\longrightarrow} 0 \\
 \implies &{|g \circ \Hat{\theta}(\Tilde{S}_n(\theta) ) - g \circ \Hat{\theta}(\bar{S}(\theta)) | \underset{n \to \infty}{\longrightarrow} 0} \, ,
 \end{split}
\end{equation}
which provides us a sufficient condition for \eqref{eq:prop_11_condition_1}.

\paragraph{Further Simplifications of the Desired Result with Successive Sufficient Conditions} We find another, simpler, sufficient condition for \eqref{eq:prop_11_condition_2} from Equation~\eqref{eq:proof_sufficient_condition}. This part is unique to our proof and absent from \cite{fort2003convergence}. It is here that we relate the formal conditions of Proposition \ref{prop:11} to the specific assumptions of our Theorem \ref{thm:main}.

We first remove the dependency on the terms $\{\theta_n\}_n$ of the EM sequence: 
\begin{equation} \label{eq:proof_intermediary_bound_supremum}
 \norm{\Tilde{S}_n(\theta_n) - \bar{S}(\theta_n) } \mathds{1}_{\theta_n \in K} \leq \sup{\theta \in K} \, \norm{\Tilde{S}_n(\theta) - \bar{S}(\theta)}.
\end{equation}

From Equations~\eqref{eq:proof_sufficient_condition}--\eqref{eq:proof_intermediary_bound_supremum} we get that:
\begin{equation*}
 \forall \text{ compact } K \subset \Theta, \quad \sup{\theta \in K} \, \norm{\Tilde{S}_n(\theta) - \bar{S}(\theta)} \underset{n \to \infty}{\longrightarrow} 0 \, ,
\end{equation*}
is a sufficient condition to have both Equations~\eqref{eq:prop_11_condition_1} and \eqref{eq:prop_11_condition_2}.

To show that the hypotheses of Theorem \ref{thm:main} imply this sufficient condition, we express it in integral form. Let $S = \brace{S_u}_{u=1, ..., q}$. We recall that $\Tilde{S}_n(\theta) = \brace{\int_z S_u(z) \Tilde{p}_{\theta, n}(z) dz}_{u=1}^q$ and $\bar{S}(\theta) = \brace{\int_z S_u(z) p_{\theta}(z) dz}_{u=1}^q$. Hence:
\begin{equation*}
 \Tilde{S}_n(\theta) - \bar{S}(\theta) = \brace{\int_z S_u(z) \parent{ \Tilde{p}_{\theta, n}(z) - p_{\theta}(z) } dz}_{u=1}^q \, .
\end{equation*}

These $q$ terms can be upper bounded by two different terms depending on the existence of the involved quantities:
\begin{equation*}
 \int_z S_u(z) \parent{ \Tilde{p}_{\theta, n}(z) - p_{\theta}(z) } dz \leq \parent{\int_z S_u(z)^2 dz}^{\frac{1}{2}} \parent{\int_z \parent{ \Tilde{p}_{\theta, n}(z) - p_{\theta}(z) }^2 dz }^{\frac{1}{2}} \, ,
\end{equation*}
and
\begin{equation*}
 \int_z S_u(z) \parent{ \Tilde{p}_{\theta, n}(z) - p_{\theta}(z) } dz \leq \parent{\int_z S_u(z)^2 p_{\theta}(z) dz }^{\frac{1}{2}} \parent{\int_z \parent{ \frac{\Tilde{p}_{\theta, n}(z)}{p_{\theta}(z)} - 1 }^2 p_{\theta}(z) dz }^{\frac{1}{2}} \, .
\end{equation*}

As a consequence, if $\int_z S_u(z)^2 dz$ exists, then it is sufficient to show have:
\begin{equation*}
 \sup{\theta \in K} \, \int_z \parent{ \Tilde{p}_{\theta, n}(z) - p_{\theta}(z) }^2 dz \underset{n \to \infty}{\longrightarrow} 0 \, ,
\end{equation*}
and if $\int_z S_u(z)^2 p_{\theta}(z) dz$ exists, then it is sufficient to show have:
\begin{equation*} 
 \sup{\theta \in K} \, \int_z \parent{ \frac{\Tilde{p}_{\theta, n}(z)}{p_{\theta}(z)} - 1 }^2 p_{\theta}(z) dz \underset{n \to \infty}{\longrightarrow} 0 \, .
\end{equation*}

Among the assumptions of Theorem \ref{thm:main} is one that states that for all compacts $K \subseteq \Theta$, one of those scenarios has to be true. Hence our sufficient condition is met.

\paragraph{Conclusions} 
With the hypotheses of Theorem \ref{thm:main}, we have 
\begin{equation*}
 \forall \text{ compact } K \subseteq \Theta, \quad \sup{\theta \in K} \, \norm{\Tilde{S}_n(\theta) - \bar{S}(\theta)} \underset{n \to \infty}{\longrightarrow} 0 \, ,
\end{equation*}
which is a sufficient condition to verify both Equation~\eqref{eq:prop_11_condition_1} and \eqref{eq:prop_11_condition_2}. With these two conditions, we can apply Proposition \ref{prop:11}.

\subsubsection{Applying Proposition \ref{prop:11}}
Since we verify all the conditions of Proposition \ref{prop:11}, we can apply its conclusions:
\begin{equation*}
 \text{With probability 1, } \underset{n \to \infty}{\lim sup}\, p_n < \infty \text{ and } \brace{\theta_n}_n \text{ compact sequence } \, ,
\end{equation*}
which is specifically the result $(i)(a)$ of Theorem \ref{thm:main}.

\subsubsection{Verifying the Conditions of Proposition \ref{prop:9}}
With Proposition \ref{prop:9}, we prove the remaining points of Theorem \ref{thm:main}: $(i)(b)$ and $(ii)$.

For the application of Proposition \ref{prop:9}:
\begin{itemize}
 \item $Cl\parent{\brace{\theta_n}_n}$ (set closure) plays the part of the compact $K$
 \item $\brace{\theta \in \Theta | T(\theta) = \theta}$ plays the part of the set $\L$
 \item $\L \cap K$ is also compact thanks to hypothesis $M3$
 \item The likelihood $g$ is the $C^0$ Lyapunov function with regards to $(T, \L)$
 \item $\brace{\theta_n}_n$ is the $K$ valued sequence (since $K$ is $Cl\parent{\brace{\theta_n}_n}$).
\end{itemize}

The last condition that remains to be shown to apply Proposition \ref{prop:9} is that:
\begin{equation*}
 \underset{n \to \infty}{\lim} | g(\theta_{n+1}) - g \circ T(\theta_n) | = 0 \, .
\end{equation*}

We have more or less already proven that, in the previous section of the Proof, with $F_n(\theta_n)$ in place of $\theta_{n+1}$. The only indices where $F_n(\theta_n) \neq \theta_{n+1} $ are when the value of the sequence $p_n$ experiences an increment of 1. We have proven with Proposition \ref{prop:11} that there is only a finite number of such increments. 
\begin{equation*}
 | g(\theta_{n+1}) - g \circ T(\theta_n) | = | g(\theta_0) - g \circ T(\theta_n) | \mathds{1}_{p_{n+1} = p_n + 1} \, + \, | g \circ F_n(\theta_n) - g \circ T(\theta_n) | \mathds{1}_{p_{n+1} = p_n} \, .
\end{equation*}

Since there is only a finite number of increments of the value of $p_n$, then $\exists N \in \N, \, \forall n \geq N, \, \mathds{1}_{p_{n+1} = p_n + 1} = 0$ and $\mathds{1}_{p_{n+1} = p_n} = 1$. In other words:
\begin{equation*}
\begin{split}
 \exists N \in \N, \, \forall n \geq N, \, \; | g(\theta_{n+1}) - g \circ T(\theta_n) | &= | g \circ F_n(\theta_n) - g \circ T(\theta_n) | \\
 \exists N \in \N, \, \forall n \geq N, \, \; | g(\theta_{n+1}) - g \circ T(\theta_n) | &= | g \circ F_n(\theta_n) - g \circ T(\theta_n) | \mathds{1}_{\theta_n \in Cl(\brace{\theta_k}_k)} \, .
\end{split}
\end{equation*}

Since $\theta_n$ is always in $Cl(\brace{\theta_k}_k)$ by definition. Additionally Proposition \ref{prop:11} tells us that $Cl(\brace{\theta_k}_k)$ is a compact. Moreover, in order to use Proposition \ref{prop:11} in the first place, we had proven that:
$${\forall \text{ compact } K \subseteq \Theta, \quad \underset{n \to \infty}{\lim} |g \circ F_n(\theta_n) - g \circ T (\theta_n)| \mathds{1}_{\theta_n \in K} = 0.}$$

We can apply this directly with $K = Cl(\brace{\theta_k}_k)$ to conclude the desired result: 
\begin{equation*}
 \underset{n \to \infty}{\lim} | g(\theta_{n+1}) - g \circ T(\theta_n) | = 0
\end{equation*}

Hence we verify all the conditions to apply Proposition \ref{prop:9}.

\subsubsection{Applying Proposition \ref{prop:9}}
Since we verify all we need, we have the conclusions of Proposition \ref{prop:9}:
\begin{itemize}
 \item $\brace{g(\theta_n)}_{n \in \N}$ converges towards a connected component of $g(\L \cap Cl(\brace{\theta_n}_n)) \subset g(\L)$
 \item If $g(\L \cap Cl(\brace{\theta_n}_n))$ has an empty interior, then $\brace{g(\theta_n)}_{n \in \N}$ converges towards a $g^* \in \R$ and $\brace{\theta_n}_n$ converges towards $\L_{g^*} \cap Cl(\brace{\theta_n}_n)$. Where $\L_{g^*} := \brace{\theta \in \L | g(\theta) = g^*}$
\end{itemize}

Both points are respectively the statements $(i)(b)$ and $(ii)$ of Theorem \ref{thm:main}.\\
Which concludes the proof of the Theorem.

\subsection{Proof of the Tempering Theorem} \label{sect:proofs_tempered}
In this Section, we prove Theorem \ref{thm:tempered_EM} of the main paper, the convergence of the tempered EM algorithm. For that, we need to show that we verify each of the hypotheses of the more general Theorem \ref{thm:main}.

We already assumed the conditions M1, M2 and M3 in the hypotheses of Theorem \ref{thm:tempered_EM}. To apply Theorem \ref{thm:main}, we need to show that when $\Tilde{p}_{\theta, n}(z) := \frac{p_{\theta}^{\frac{1}{T_n}}(z)}{\int_{z'} p_{\theta}^{\frac{1}{T_n}}(z') dz'}$, then $\forall \text{ compact } K \subseteq \Theta$, one of the two following configurations holds:
\begin{equation*} 
 \int_z S(z)^2 dz < \infty \text{ and }\sup{\theta \in K} \, \int_z \parent{ \Tilde{p}_{\theta, n}(z) - p_{\theta}(z) }^2 dz \underset{n \to \infty}{\longrightarrow} 0 \, , 
\end{equation*}
or
\begin{equation*} 
 \sup{\theta \in K} \, \int_z S(z)^2 p_{\theta}(z) dz < \infty \text{ and } \sup{\theta \in K} \, \int_z \parent{ \frac{\Tilde{p}_{\theta, n}(z)}{p_{\theta}(z)} - 1 }^2 p_{\theta}(z) dz \underset{n \to \infty}{\longrightarrow} 0 \, .
\end{equation*}

Since we have assumed:
\begin{equation*}
 \forall \text{ compact } K \in \Theta, \; \forall \alpha \in \overline{\mathcal{B}}(1, \epsilon), \forall u, \quad \sup{\theta \in K}\, \int_z S^2_u(z) p_{\theta}^{\alpha}(z) dz < \infty \, , 
\end{equation*}
then we already verify the first half of the second configuration for all the compacts $K$. Hence it is sufficient to prove that:
\begin{equation} \label{eq:thm_2_proof_intermediary_result}
 \forall \text{ compact } K \in \Theta, \; \sup{\theta \in K} \, \int_z \parent{ \frac{\Tilde{p}_{\theta, n}(z)}{p_{\theta}(z)} - 1 }^2 p_{\theta}(z) dz \underset{n \to \infty}{\longrightarrow} 0 \, , 
\end{equation}
to have the desired result. The rest of the proof is dedicated to this goal.

\subsubsection{Taylor Development}
We use the Taylor's formula of the first order with the mean-value form of the reminder. For a derivable function $f$:
\begin{equation} \label{eq:proof_taylor}
 f(x) = f(0) + f'(a) x , \quad a \in \brack{0, x},
\end{equation}
where the interval $ \brack{0, x}$ has a flexible meaning since $x$ could be negative.

We apply it to:
\begin{equation*}
 f(x) = e^x, \quad f'(x) = e^x, \quad f(x) = 1 + x e^a, \; a \in \brack{0, x} \, , 
\end{equation*}
and:
\begin{equation*}
 f(x) = \frac{1}{1+x}, \quad f'(x) = -\frac{1}{(1+x)^2}, \quad f(x) = 1 - \frac{x}{(1+a)^2}, \; a \in \brack{0, x} \, .
\end{equation*}

To make the upcoming calculation more readable, we momentarily replace $p_{\theta}(z)$ by simply $p$ and $T_n$ by T.
\begin{equation*}
 \begin{split}
 p^{\frac{1}{T}} &= p\parent{ p^{\frac{1}{T}-1} }\\
 &= p e^{(\frac{1}{T}-1)\ln p}\\
 &= p + \parent{\frac{1}{T}-1} p \,\ln p \, e^{a}, \quad a \in \brack{0, \parent{\frac{1}{T}-1} \ln p}\, ,
 \end{split}
\end{equation*}
where $a = a(z, \theta, T_n)$ since it depends on the value of $p_{\theta}(z)$ and $T_n$. Provided that the following quantities are defined, we have:
\begin{equation*}
 \int_z p^{\frac{1}{T}} = 1 + \parent{\frac{1}{T}-1} \int_z p \,\ln p \, e^{a}\, ,
\end{equation*}
Hence:
\begin{equation*}
 \frac{1}{\int_z p^{\frac{1}{T}}} = 1 - \parent{\frac{1}{T}-1} \frac{\int_z p \,\ln p \, e^{a}}{\parent{1 + b}^2}, \quad b \in \brack{0, \parent{\frac{1}{T}-1} \int_z p \,\ln p \, e^{a}}\, ,
\end{equation*}
where $b = b(\theta, T_n)$ since it depends on the value of $T_n$ the integral over $z$ of a function of $z$ and $\theta$.

In the end, we have:

\begin{equation} \label{eq:proof_full_taylor_ld}
 \frac{p^{\frac{1}{T}}}{\int_z p^{\frac{1}{T}}} = p + \parent{\frac{1}{T}-1} p \,\ln p \, e^{a} \parent{1 - \parent{\frac{1}{T}-1} \frac{\int_z p \,\ln p \, e^{a}}{\parent{1 + b}^2}} - \parent{\frac{1}{T}-1} p \frac{\int_z p \,\ln p \, e^{a}}{\parent{1 + b}^2} \, .
\end{equation}

Since for any real numbers $(x+y)^2 \leq 2 (x^2 + y^2)$, then:
\begin{equation*}
\begin{split}
 \parent{ \frac{p^{\frac{1}{T}}}{\int_z p^{\frac{1}{T}}} \!-\! p }^2 \!&\leq 2\parent{\frac{1}{T}-1}^2 \!p^2 \parent{ \parent{\ln p \, e^{a}}^2 \parent{1 \!-\! \parent{\frac{1}{T}-1} \frac{\int_z p \,\ln p \, e^{a}}{\parent{1 + b}^2}}^2 \!+\! \parent{\frac{\int_z p \,\ln p \, e^{a}}{\parent{1 + b}^2}}^2 } \\
 &= 2\parent{\frac{1}{T}-1}^2 p^2 \parent{ \parent{\ln p \, e^{a}}^2 A + B } \, .
\end{split}
\end{equation*}
where $A = A(\theta, T_n)$ and $B = B(\theta, T_n)$. So far the only condition that has to be verified for all the involved quantities to be defined is that $\int_z p \,\ln p \, e^{a}$ exists. With this Taylor development on hand, we state, prove and apply two lemmas which allow us to get \eqref{eq:thm_2_proof_intermediary_result} and conclude the proof of the theorem.

\subsubsection{Two Intermediary Lemmas}
The two following lemmas provides every result we need to finish the proof.
\begin{Lemma} \label{lem:link_S_log} With 
 \begin{equation*}
 p_{\theta}(z) = exp\parent{\psi(\theta) + \dotprod{S(z), \phi(\theta}}\, ,
 \end{equation*}
 then
 \begin{equation*}
 \int_z p_{\theta}^{\alpha}(z) \, \ln^2 \, p_{\theta}(z) dz \leq 2 \psi(\theta)^2 \int_z p_{\theta}^{\alpha}(z) dz + 2 \norm{\phi(\theta)}^2 . \sum_u \int_z S^2_u(z) p_{\theta}^{\alpha}(z) \, .
 \end{equation*}
 and
 \begin{equation*}
 \int_z p_{\theta}^{\alpha}(z) \, \det{\ln p_{\theta}(z)} dz \leq \det{\psi(\theta)} \int_z p_{\theta}^{\alpha}(z) dz + \norm{\phi(\theta)} . \parent{ \sum_u \int_z S^2_u(z) p_{\theta}^{\alpha}(z) \int_z p_{\theta}^{\alpha}(z) }^{\frac{1}{2}} \, .
 \end{equation*}
 
\end{Lemma}

\begin{proof}
 For the first inequality, using the fact that $(a+b)^2 \leq 2(a^2+b^2)$, we have:
 \begin{equation*}
 \int_z p_{\theta}^{\alpha}(z) \, \ln^2 \, p_{\theta}(z) dz \leq 2 \psi(\theta)^2 \int_z p_{\theta}^{\alpha}(z) dz + 2 \int_z p_{\theta}^{\alpha}(z) \dotprod{S(z), \phi(\theta)}^2 \, ,
 \end{equation*}
 We use Cauchy-Schwartz:
 \begin{equation*}
 \dotprod{S(z), \phi(\theta)}^2 \leq \norm{\phi}^2 \norm{S(z)}^2 = \norm{\phi}^2 \sum_u S_u(z)^2\, ,
 \end{equation*}
 to get the desired result:
 \begin{equation*}
 \int_z p_{\theta}^{\alpha}(z) \, \ln^2 \, p_{\theta}(z) dz \leq 2 \psi(\theta)^2 \int_z p_{\theta}^{\alpha}(z) dz + 2 \norm{\phi(\theta)}^2 . \sum_u \int_z S^2_u(z) p_{\theta}^{\alpha}(z) \, .
 \end{equation*}
 For the second inequality, we start with Cauchy-Schwartz on $\dotprod{\int_z S(z) p_{\theta}^{\alpha}(z), \phi(\theta)}$:
 \begin{equation*}
 \int_z p_{\theta}^{\alpha}(z) \, \det{\ln p_{\theta}(z)} dz \leq \det{\psi(\theta)} \int_z p_{\theta}^{\alpha}(z) dz + \norm{\phi(\theta)} . \norm{\int_z S(z) p_{\theta}^{\alpha}(z)} \, .
 \end{equation*}
 Moreover, since:
 \begin{equation*}
 \int_z S_u(z) p_{\theta}^{\alpha}(z) dz \leq \parent{\int_z S^2_u(z) p_{\theta}^{\alpha}(z) dz }^{\frac{1}{2}} \parent{\int_z p_{\theta}^{\alpha}(z) dz }^{\frac{1}{2}}\, ,
 \end{equation*}
 then
 \begin{equation*}
 \int_z p_{\theta}^{\alpha}(z) \, \det{\ln p_{\theta}(z)} dz \leq \det{\psi(\theta)} \int_z p_{\theta}^{\alpha}(z) dz + \norm{\phi(\theta)} . \parent{ \sum_u \int_z S^2_u(z) p_{\theta}^{\alpha}(z) \int_z p_{\theta}^{\alpha}(z) }^{\frac{1}{2}} \, .
 \end{equation*}
\end{proof}

\begin{Lemma} \label{lem:scenario_2} With $K$ compact and $\epsilon \in \R_+^*$,
 \begin{equation*}
 p_{\theta}(z) = exp\parent{\psi(\theta) + \dotprod{S(z), \phi(\theta}},
 \end{equation*}
 and
 \begin{equation*}
 \Tilde{p}_{\theta, n}(z) := \frac{p_{\theta}^{\frac{1}{T_n}}(z)}{\int_{z'} p_{\theta}^{\frac{1}{T_n}}(z') dz'},
 \end{equation*}
 if 
 \begin{enumerate}[(i)]
 \item $T_n \in \R^*_+ \underset{n \to \infty}{\longrightarrow} 1$ \, ,
 \item $ \sup{\theta \in K} \,\psi(\theta) < \infty $ \, , 
 \item $ \sup{\theta \in K}\, \norm{\phi(\theta)} < \infty $ \, , 
 
 \item $\forall \alpha \in \overline{\mathcal{B}}(1, \epsilon), \quad \sup{\theta \in K}\, \int_z p_{\theta}^{\alpha}(z) dz < \infty$ \, , 
 \item $\forall \alpha \in \overline{\mathcal{B}}(1, \epsilon), \; \forall u, \quad \sup{\theta \in K}\, \int_z S^2_u(z) p_{\theta}^{\alpha}(z) dz < \infty$ \, .
 \end{enumerate}
 then 
 $$ \sup{\theta \in K} \, \int_z \parent{ \frac{\Tilde{p}_{\theta, n}(z)}{p_{\theta}(z)} - 1 }^2 p_{\theta}(z) dz \underset{n \to \infty}{\longrightarrow} 0 \, . $$
\end{Lemma}

\begin{proof}
Provided that the following integrals exist, we have, thanks to the Taylor development:
\begin{equation} \label{eq:proof_lemma_scenario_2_big_inequality}
\begin{split}
 \int_z \frac{1}{p} \parent{ \frac{p^{\frac{1}{T}}}{\int_z p^{\frac{1}{T}}} - p }^2 &\leq 2 \int_z \parent{\frac{1}{T}-1}^2 p \parent{ \parent{\ln p \, e^{a}}^2 A + B } \\
 &= 2 \parent{\frac{1}{T}-1}^2 A \int_z p e^{2a} \ln^2 \, p \, + 2 \parent{\frac{1}{T}-1}^2 B \, .
\end{split}
\end{equation}
In this proof, we find finite upper bounds independent of $\theta$ and $T_n$ for $A(\theta, T_n)$, $B(\theta, T_n)$ and $\int_z p e^{2a} \ln^2 \, p$, then - since $\parent{\frac{1}{T_n}-1} \longrightarrow 0$ - we have the desired result.\\
\\
We start by studying $A(\theta, T_n) = \parent{1 - \parent{\frac{1}{T}-1} \frac{\int_z p \,\ln p \, e^{a}}{\parent{1 + b}^2}}^2$. The first term of interest here is $\int_z p \,\ln p \, e^{a}$. We have:
\begin{equation*}
\begin{split}
 a &\in \brack{0, \parent{\frac{1}{T}-1} \ln p}\, , \\
 e^a &\in \brack{1, p^{\frac{1}{T}-1} } \, , \\
 p \,\ln p \, e^{a} &\in \brack{p \, \ln p, p^{\frac{1}{T}} \ln p} \, .
\end{split}
\end{equation*}
where we recall that the interval is to be taken in a flexible sense, since we do not now a priory which bound is the largest and which is the smallest. What we have without doubt though is:
\begin{equation*}
 \det{p \,\ln p \, e^{a}} \leq max\parent{ \det{p \, \ln p} , \det{p^{\frac{1}{T}} \ln p}} \, .
\end{equation*}
We find an upper bound on both those term. Let $\alpha \in \overline{\mathcal{B}}(1, \epsilon)$, the second result of Lemma~\ref{lem:link_S_log} provides us:
\begin{equation*}
 \int_z p_{\theta}^{\alpha}(z) \, \det{\ln p_{\theta}(z)} dz \leq \det{\psi(\theta)} \int_z p_{\theta}^{\alpha}(z) dz + \norm{\phi(\theta)} . \parent{ \sum_u \int_z S^2_u(z) p_{\theta}^{\alpha}(z) \int_z p_{\theta}^{\alpha}(z) }^{\frac{1}{2}} \, .
\end{equation*}
Thanks to the hypotheses (ii), (iii), (iv) and (v), we have:
\begin{equation*}
\begin{split}
 \int_z p_{\theta}^{\alpha}(z) \, \det{\ln p_{\theta}(z)} dz &\leq \sup{\theta \in K}\, \det{\psi(\theta)} . \sup{\theta \in K}\, \int_z p_{\theta}^{\alpha}(z) dz \\
 &\hspace{0.5cm}+ \sup{\theta \in K}\, \norm{\phi(\theta)} . \sum_u \parent{ \sup{\theta \in K}\, \int_z S^2_u(z) p_{\theta}^{\alpha}(z)}^{\frac{1}{2}}. \parent{ \sup{\theta \in K}\, \int_z p_{\theta}^{\alpha}(z) }^{\frac{1}{2}}\\
 &=: C(\alpha)\\
 &<\infty\, .
\end{split}
\end{equation*}

The upper bound $C(\alpha)$ in the previous inequality is independent of $\theta$ and $z$ but still dependant of the exponent $\alpha$. However, since $\overline{\mathcal{B}}(1, \epsilon)$ is closed ball, hypotheses (iv) and (v) can be rephrased as:
\begin{equation*}
 \begin{split}
 (iv)& \; \sup{\alpha \in \overline{\mathcal{B}}(1, \epsilon)}\, \sup{\theta \in K}\, \int_z p_{\theta}^{\alpha}(z) dz < \infty \, ,\\
 (v)& \; \forall u, \; \sup{\alpha \in \overline{\mathcal{B}}(1, \epsilon)}\, \sup{\theta \in K}\, \int_z S^2_u(z) p_{\theta}^{\alpha}(z) dz < \infty\, .
 \end{split}
\end{equation*}

Hence we can actually take the supremum over $\alpha$ in the right term of the inequation as well. We have: $\int_z p_{\theta}^{\alpha}(z) \, \det{\ln p_{\theta}(z)} dz$
\begin{equation*}
\begin{split}
 &\leq \sup{\theta \in K}\, \det{\psi(\theta)} . \sup{\alpha \in \overline{\mathcal{B}}(1, \epsilon)}\, \sup{\theta \in K}\, \int_z p_{\theta}^{\alpha}(z) dz \\
 &\hspace{0.5cm}+ \sup{\theta \in K}\, \norm{\phi(\theta)} . \sum_u \parent{ \sup{\alpha \in \overline{\mathcal{B}}(1, \epsilon)}\, \sup{\theta \in K}\, \int_z S^2_u(z) p_{\theta}^{\alpha}(z)}^{\frac{1}{2}}. \parent{\sup{\alpha \in \overline{\mathcal{B}}(1, \epsilon)}\, \sup{\theta \in K}\, \int_z p_{\theta}^{\alpha}(z) }^{\frac{1}{2}}\\
 & =:C'\\
 &<\infty\, .
\end{split}
\end{equation*}

This new upper bound $C'$ is independent of $\alpha$.

Since $T_n \mapsto 1$, then $ \exists N \in \N, \, \forall n \geq N, \, \frac{1}{T_n} \in \overline{\mathcal{B}}(1, \epsilon) $. Hence for $n \geq N$, we can apply the previous inequation to either $\alpha = 1$ or $\alpha = \frac{1}{T_n}$, which provides us that $\int_z p_{\theta}(z) \, \det{\ln p_{\theta}(z)}$, $\int_z p_{\theta}^{\frac{1}{T_n}}(z) \, \det{\ln p_{\theta}(z)}$ and their supremum in $\theta$ are all finite, all of them upper bounded by $C'$. \\

In the end, when $n\geq N$, we have the control $\sup{\theta \in K} \, \det{\int_z p \,\ln p \, e^{a}} < C'$. \\

The next term to control is $\frac{1}{\parent{1 + b}^2}$.

Since $b \in \brack{0, \parent{\frac{1}{T}-1} \int_z p \,\ln p \, e^{a}},$ then $|b| \leq \parent{\frac{1}{T}-1} \sup{\theta \in K} \, \int_z p \,\ln p \, e^{a}$. We already established that for all $n \geq N$, $\sup{\theta \in K} \, \det{\int_z p \,\ln p \, e^{a}} \leq C' < \infty$, hence $\sup{\theta \in K}\, |b(\theta, T_n)| \underset{T_n \longrightarrow 1}{\longrightarrow}0$. In particular, $\exists N' \in \N, \forall n \geq N', \forall \theta \in K$ we have $|b(\theta, T_n)| \leq \frac{1}{2}$. In that case:
\begin{equation*}
 \begin{split}
 \parent{1+b}^2 > \parent{1-|b|}^2 \geq \frac{1}{4}\\
 \frac{1}{\parent{1+b}^2} < \frac{1}{\parent{1-|b|}^2}\leq 4 \, .
 \end{split}
\end{equation*}

In the end, when $n \geq max(N, N') $, for any $\theta \in K$:
\begin{equation*}
 \begin{split}
 A(\theta, T_n) &\leq 2 + 2 \parent{\frac{1}{T_n} - 1}^2 \parent{\frac{\int_z p \,\ln p \, e^{a}}{\parent{1 + b}^2}}^2 \\
 &\leq 2 + 32 \parent{\frac{1}{T_n} - 1}^2 \parent{\sup{\theta \in K} \, \int_z p \,\ln p \, e^{a}}^2 \\
 &\leq 2 + 32 \parent{\frac{1}{T_n} - 1}^2 C'^2\\
 & \leq 2 + 32 \epsilon^2 C'^2\\
 &=: C_1\, .
 \end{split}
\end{equation*}

This upper bound does not depend en $\theta$ anymore and the part in $T_n$ simply converges towards 0 when $T_n \longrightarrow 1$.

Treating the term $B(\theta, T_n) = \parent{\frac{\int_z p \,\ln p \, e^{a}}{\parent{1 + b}^2}}^2 \leq 16 \parent{\sup{\theta \in K} \, \int_z p \,\ln p \, e^{a}}^2 \leq 16 C'^2 =: C_2$ is immediate after having dealt with $A(\theta, T_n)$.

We now treat the term $\int_z p \, e^{2a} \ln^2 \, p$ in the exact same fashion as we did $A(\theta, T_n)$: 
\begin{equation*}
 \begin{split}
 p \,\ln p \, e^{a} &\in \brack{p \, \ln p, p^{\frac{1}{T}} \ln p}\\
 \implies p \parent{\ln p \, e^{a}}^2 &\in \brack{p \, \ln^2 \, p, p^{\frac{2}{T}-1} \ln^2 \, p}\\
 \implies p \parent{\ln p \, e^{a}}^2 &\leq max(p \, \ln^2 \, p, p^{\frac{2}{T}-1} \ln^2 \, p) \, .
 \end{split}
\end{equation*}

We control those two terms as previously. First we apply Lemma \ref{lem:link_S_log} (its first result this time) with $\alpha \in \overline{\mathcal{B}}(1, \epsilon)$. 
\begin{equation*}
 \int_z p_{\theta}^{\alpha}(z) \, \ln^2 \, p_{\theta}(z) dz \leq 2 \psi(\theta)^2 \int_z p_{\theta}^{\alpha}(z) dz + 2 \norm{\phi(\theta)}^2 . \sum_u \int_z S^2_u(z) p_{\theta}^{\alpha}(z) \, .
\end{equation*}

Thanks to the hypotheses (ii), (iii), (iv) and (v), we can once again take the supremum of the bound over $\theta \in K$, then over $\alpha \in \overline{\mathcal{B}}(1, \epsilon)$ and conserve finite quantities:

\begin{equation*}
\begin{split}
 \int_z p_{\theta}^{\alpha}(z) \, \ln^2 \, p_{\theta}(z) dz &\leq 2 \sup{\theta \in K}\, \psi(\theta)^2 . \sup{\alpha \in \overline{\mathcal{B}}(1, \epsilon)}\, \sup{\theta \in K}\, \int_z p_{\theta}^{\alpha}(z) dz \\
 & \hspace{0.2cm} + 2 \sup{\theta \in K}\, \norm{\phi(\theta)}^2 . \sum_u \sup{\alpha \in \overline{\mathcal{B}}(1, \epsilon)}\, \sup{\theta \in K}\, \int_z S^2_u(z) p_{\theta}^{\alpha}(z) \\
 &=: C_3 \\
 &< \infty \, .
\end{split}
\end{equation*}

The previous result is true for $\alpha=1$, and since once again $\exists N'', \, \forall n \geq N'', \frac{2}{T_n}-1 \in \overline{\mathcal{B}}(1, \epsilon) \cap \R^*_+ $, it is also true for $\alpha = \frac{2}{T_n}-1$ when $n$ is large enough. $C_3$ is independent of $z$, $\theta$ and $T_n$.

In the end $\forall n\geq N'', \; \int_z p \, e^{2a} \ln^2 \, p \leq C_3 < \infty $.

We replace the three terms $A(\theta, T_n)$, $B(\theta, T_n)$ and $\int_z p \, e^{2a} \ln^2 \, p$ by their upper bounds in the inequality \eqref{eq:proof_lemma_scenario_2_big_inequality}. When $n \geq max(N, N', N'')$:
\begin{equation*}
 \int_z \frac{1}{p} \parent{ \frac{p^{\frac{1}{T}}}{\int_z p^{\frac{1}{T}}} - p }^2 \leq 2 \parent{\frac{1}{T_n}-1}^2 C_1 C_3 + 2 \parent{\frac{1}{T_n}-1}^2 C_2\, .
\end{equation*}

Which converges towards 0 when $T_n \longrightarrow 1$, i.e. when $n\longrightarrow \infty$. This concludes the proof of the lemma.

\end{proof}

\subsubsection{Verifying the Conditions of Lemma \ref{lem:scenario_2} }
Now that the lemmas are proven, all that remains is to apply Lemma \ref{lem:scenario_2}.\\
\paragraph{\textit{(i)}} We have $T_n \in \R^*_+ \underset{n \to \infty}{\longrightarrow} 1$ by hypothesis.
\paragraph{\textit{(ii)} and \textit{(iii)}} $\sup{\theta \in K} \,\psi(\theta) < \infty $ and $ \sup{\theta \in K}\, \norm{\phi(\theta)} < \infty $ are implied by the fact that $\psi(\theta) = \psi'(\theta) - \ln g(\theta)$ and $\phi(\theta)$ are continuous
\paragraph{\textit{(iv)} and \textit{(v)}} are also hypotheses of the theorem.

Hence we can apply Lemma \ref{lem:scenario_2}. This means that:
$$ \sup{\theta \in K} \, \int_z \parent{ \frac{\Tilde{p}_{\theta, n}(z)}{p_{\theta}(z)} - 1 }^2 p_{\theta}(z) dz \underset{n \to \infty}{\longrightarrow} 0 . $$
with this last condition verified, we can apply Theorem \ref{thm:main}, which concludes the proof.

\end{document}